\documentclass[a4paper]{article}
\usepackage{amssymb}
\usepackage{amsmath}
\usepackage{amsthm}
\usepackage{authblk}
\usepackage[USenglish]{babel}
\usepackage{esint}
\usepackage{graphicx}
\usepackage[colorlinks=true]{hyperref}
\usepackage[latin1]{inputenc}
\newtheorem{theorem}{Theorem}[section]
\newtheorem{lemma}[theorem]{Lemma}
\newtheorem{proposition}[theorem]{Proposition}
\newtheorem{corollary}[theorem]{Corollary}
\newtheorem{remark}[theorem]{Remark}
\newtheorem{definition}[theorem]{Definition}

\newcommand{\N}{\mathbb N}
\newcommand{\Z}{\mathbb Z}

\newcommand{\R}{\mathbb R}

\newcommand{\B}{\mathbb B}

\newcommand{\mcal}{\mathcal}

\newcommand{\mrm}{\mathrm}

\renewcommand{\a}{\alpha}

\newcommand{\g}{\gamma}
\newcommand{\G}{\Gamma}
\renewcommand{\d}{\delta}
\newcommand{\D}{\Delta}
\newcommand{\e}{\varepsilon}
\newcommand{\z}{\zeta}
\renewcommand{\t}{\theta}

\newcommand{\la}{\lambda}
\newcommand{\La}{\Lambda}
\newcommand{\s}{\sigma}

\newcommand{\Si}{\Sigma}

\newcommand{\ph}{\varphi}

\renewcommand{\O}{\Omega}
\newcommand{\wt}{\widetilde}

\newcommand{\ol}{\overline}

\newcommand{\ub}{\underbrace}
\newcommand{\fr}{\frac}

\newcommand{\n}{\nabla}
\newcommand{\fa}{\forall}
\newcommand{\ex}{\exists}
\newcommand{\es}{\emptyset}
\newcommand{\wk}{\rightharpoonup}
\newcommand{\inc}{\hookrightarrow}

\newcommand{\us}{\underset}
\newcommand{\sr}{\stackrel}

\newcommand{\To}{\Rightarrow}

\newcommand{\lto}{\longrightarrow}

\newcommand{\sm}{\setminus}
\renewcommand{\Cup}{\bigcup}
\renewcommand{\Cap}{\bigcap}
\newcommand{\sub}{\subset}
\newcommand{\Sub}{\Subset}

\newcommand{\nin}{\not\in}
\newcommand{\eq}{\equiv}

\newcommand{\pl}{\oplus}

\newcommand{\x}{\times}
\renewcommand{\c}{\circ}
\newcommand{\cd}{\cdot}
\newcommand{\ds}{\dots}

\newcommand{\tx}{\text}
\newcommand{\q}{\quad}
\renewcommand{\l}{\left}
\renewcommand{\r}{\right}

\newcommand{\bthm}{\begin{theorem}}
\newcommand{\ethm}{\end{theorem}}
\newcommand{\blem}{\begin{lemma}}
\newcommand{\elem}{\end{lemma}}
\newcommand{\bprop}{\begin{proposition}}
\newcommand{\eprop}{\end{proposition}}
\newcommand{\bcor}{\begin{corollary}}
\newcommand{\ecor}{\end{corollary}}
\newcommand{\bdefi}{\begin{definition}}
\newcommand{\edefi}{\end{definition}}
\newcommand{\bpf}{\begin{proof}}
\newcommand{\epf}{\end{proof}}
\newcommand{\bl}{\begin{array}{l}}
\newcommand{\bll}{\begin{array}{ll}}
\newcommand{\barr}{\begin{array}}
\newcommand{\earr}{\end{array}}
\newcommand{\bite}{\begin{itemize}}
\newcommand{\eite}{\end{itemize}}
\newcommand{\bequ}{\begin{equation}}
\newcommand{\eequ}{\end{equation}}
\newcommand{\beqa}{\begin{eqnarray}}
\newcommand{\eeqa}{\end{eqnarray}}
\newcommand{\beqy}{\begin{eqnarray*}}
\newcommand{\eeqy}{\end{eqnarray*}}
\newcommand{\bin}[2]{\left(\genfrac{}{}{0pt}{}{#1}{#2}\right)}
\newcommand{\qm}[1]{``#1''}

\begin{document}

\everymath{\displaystyle}

\title{$B_2$ and $G_2$ Toda systems on compact surfaces: a variational approach}
\author{Luca Battaglia\thanks{Sapienza Universit\`a di Roma, Dipartimento di Matematica, Piazzale Aldo Moro $5$, $00185$ Roma - battaglia@mat.uniroma1.it}}
\date{}

\maketitle\

\begin{abstract}
\noindent We consider the $B_2$ and $G_2$ Toda systems on a compact surface $(\Si,g)$
$$\l\{\bl-\D u_1=a_{11}\rho_1\l(\fr{h_1e^{u_1}}{\int_\Si h_1e^{u_1}\mrm dV_g}-1\r)+a_{12}\rho_2\l(\fr{h_2e^{u_2}}{\int_\Si h_2e^{u_2}\mrm dV_g}-1\r)\\-\D u_2=a_{21}\rho_1\l(\fr{h_1e^{u_1}}{\int_\Si h_1e^{u_1}\mrm dV_g}-1\r)+a_{22}\rho_2\l(\fr{h_2e^{u_2}}{\int_\Si h_2e^{u_2}\mrm dV_g}-1\r)\earr\r.,$$
where $A=(a_{ij})=\l(\barr{cc}2&-1\\-2&2\earr\r)$ or $\l(\barr{cc}2&-1\\-3&2\earr\r)$ and $h_i\in C^\infty_{>0}(\Si),\,\rho_i\in\R_{>0}$ are given.\\
We attack the problem using variational techniques, following the previous work \cite{bjmr} concerning the $A_2$ Toda system, namely the case $A=\l(\barr{cc}2&-1\\-1&2\earr\r)$. We get existence and multiplicity of solutions as long as $\chi(\Si)\le0$ and $\rho_1,\rho_2\nin4\pi\N$.\\
We also extend some of the results to the case of general systems.\\
\end{abstract}\

\section{Introduction}\

Let $(\Si,g)$ be a closed Riemann surface with surface area equal to $1$. The $B_2$ and $G_2$ Toda systems are respectively the following systems of PDEs on $\Si$:

\bequ\label{b2toda}
\l\{\bl-\D u_1=2\rho_1\l(\fr{h_1e^{u_1}}{\int_\Si h_1e^{u_1}\mrm dV_g}-1\r)-\rho_2\l(\fr{h_2e^{u_2}}{\int_\Si h_2e^{u_2}\mrm dV_g}-1\r)\\-\D u_2=2\rho_2\l(\fr{h_2e^{u_2}}{\int_\Si h_2e^{u_2}\mrm dV_g}-1\r)-2\rho_1\l(\fr{h_1e^{u_1}}{\int_\Si h_1e^{u_1}\mrm dV_g}-1\r)\earr\r.;
\eequ

\bequ\label{g2toda}
\l\{\bl-\D u_1=2\rho_1\l(\fr{h_1e^{u_1}}{\int_\Si h_1e^{u_1}\mrm dV_g}-1\r)-\rho_2\l(\fr{h_2e^{u_2}}{\int_\Si h_2e^{u_2}\mrm dV_g}-1\r)\\-\D u_2=2\rho_2\l(\fr{h_2e^{u_2}}{\int_\Si h_2e^{u_2}\mrm dV_g}-1\r)-3\rho_1\l(\fr{h_1e^{u_1}}{\int_\Si h_1e^{u_1}\mrm dV_g}-1\r)\earr\r..
\eequ

Here, $-\D=-\D_g$ is the Laplace-Beltrami operator, $\rho_1,\rho_2$ are positive parameters and $h_1,h_2$ are positive smooth functions on $\Si$.\\

Such systems are particularly interesting because their matrices of coefficients
$$B_2=\l(\barr{cc}2&-1\\-2&2\earr\r)\q\q\q G_2=\l(\barr{cc}2&-1\\-3&2\earr\r)$$
are the Cartan matrices of the special orthonormal group $SO(5)$ and of the symplectic group $Sp(4)$, respectively. Together with $A_2=\l(\barr{cc}2&-1\\-1&2\earr\r)$, corresponding to $SU(3)$, these are the only $2$-dimensional Cartan matrices.\\
Such Toda systems have important applications in both algebraic geometry and mathematical physics. In geometry, they appear in the study of complex holomorphic curves (see \cite{cal,dunne,bw,dol}); in physics, they arise in non-Abelian gauge field theory (see \cite{lee,kaol,dunne,dun2,dun3}).\\

Both \eqref{b2toda} and \eqref{g2toda} are variational problems. In fact, solutions are respectively critical points of the following energy functionals:

\bequ\label{jb2}
J_{B_2,\rho}(u):=\int_\Si Q_{B_2}(u)\mrm dV_g-\rho_1\l(\log\int_\Si h_1e^{u_1}\mrm dV_g-\int_\Si u_1\mrm dV_g\r)-\fr{\rho_2}2\l(\log\int_\Si h_2e^{u_2}\mrm dV_g-\int_\Si u_2\mrm dV_g\r);
\eequ

\bequ\label{jg2}
J_{G_2,\rho}(u):=\int_\Si Q_{G_2}(u)\mrm dV_g-\rho_1\l(\log\int_\Si h_1e^{u_1}\mrm dV_g-\int_\Si u_1\mrm dV_g\r)-\fr{\rho_2}3\l(\log\int_\Si h_2e^{u_2}\mrm dV_g-\int_\Si u_2\mrm dV_g\r).
\eequ

Here, $Q_{B_2}$ and $Q_{G_2}$ are defined by
$$Q_{B_2}(u)=\fr{|\n u_1|^2}2+\fr{\n u_1\cd\n u_2}2+\fr{|\n u_2|^2}4,\q\q\q Q_{G_2}(u)=|\n u_1|^2+\n u_1\cd\n u_2+\fr{|\n u_2|^2}3,$$
$\n=\n_g$ is the gradient given by the metric $g$ and $\cd$ is the Riemannian scalar product.\\

A first tool to study the variational properties of systems \eqref{b2toda} and \eqref{g2toda} is given by the Moser-Trudinger type inequality proved in \cite{bat2}.\\
It basically says that both the functionals $J_{B_2,\rho}$ and $J_{G_2,\rho}$ are bounded from below on $H^1(\Si)^2$ if and only if $\rho_1,\rho_2\le4\pi$ and that they are coercive if and only if both strict inequalities holds. As a consequence, in the latter case \eqref{b2toda} and \eqref{g2toda} admit energy-minimizing solutions.\\

The purpose of this paper is to prove the existence of solutions when there are no hopes of finding global minimizers, namely when $\rho_1$ and/or $\rho_2$ exceed $4\pi$. We will look for min-max solution by using a Morse-theoretical argument based on the topology of very low sub-levels of the energy functional.\\
To apply such arguments we will need some compactness conditions on the solutions of \eqref{b2toda}, \eqref{g2toda}, which were recently proved in \cite{lz}. This result, combined with a standard monotonicity trick from \cite{luc,str}, allows to apply such min-max methods for the problem \eqref{b2toda} as long as neither $\rho_1$ nor $\rho_2$ are integer multiples of $4\pi$; the same holds true for \eqref{g2toda} under assuming an extra upper bound on both parameters.\\

Min-max methods have been used several times to get existence results for the $A_2$ Toda system
\bequ\label{a2toda}
\l\{\bl-\D u_1=2\rho_1\l(\fr{h_1e^{u_1}}{\int_\Si h_1e^{u_1}\mrm dV_g}-1\r)-\rho_2\l(\fr{h_2e^{u_2}}{\int_\Si h_2e^{u_2}\mrm dV_g}-1\r)\\-\D u_2=2\rho_2\l(\fr{h_2e^{u_2}}{\int_\Si h_2e^{u_2}\mrm dV_g}-1\r)-\rho_1\l(\fr{h_1e^{u_1}}{\int_\Si h_1e^{u_1}\mrm dV_g}-1\r)\earr\r.
\eequ
and for the well-known scalar Liouville equation
\bequ\label{liouscal}
-\D u=2\rho\l(\fr{he^u}{\int_\Si he^u\mrm dV_g}-1\r),
\eequ
which can be found back in either \eqref{b2toda}, \eqref{g2toda} or \eqref{a2toda} by setting $\rho_2:=0$.\\
A general existence result have been given in \cite{dm,dja} for \eqref{liouscal}, as well as some existence results for \eqref{a2toda} under an upper bound on one or both the $\rho_i$ (see \cite{mn,mr13,jkm}). Moreover, in \cite{bjmr} existence of solutions is shown under no assumptions on $\rho$ but rather on the topology of $\Si$.\\

In this paper we show for the first time the existence of min-max solutions for the $B_2$ and $G_2$ Toda system, by extending the general result from \cite{bjmr}.\\
In other words, we show existence of solution under only assuming $\chi(\Si)\le0$. No condition is required on $\rho$ other than satisfying the necessary compactness assumptions.\\

The assumptions on $\Si$, which is satisfied as long as $\Si$ is not homeomorphic to a sphere nor to a projective plane, enables to build two surjective retractions $\Pi_i:\Si\to\g_i$ on disjointed simple closed curves. Such retractions simplify the analysis of the energy sub-levels, essentially because the interaction of the two components is not seen, being $\g_1\cap\g_2=\es$.\\
Although retracting causes a loss of topological information, we can still prove that very low sub-levels are not contractible, which is what is needed to obtain existence of solutions.\\
Roughly speaking, the interaction between $u_1$ and $u_2$ is ruled by the matrix entries $a_{12},a_{21}$. Therefore, since this argument does not take account such an interaction, it can be extended from the case of $A_2$ matrix to $B_2,G_2$ matrices, which differ from $A_2$ just by the coefficient $a_{21}$.\\
In general, the interaction between components should be quite hard to handle. Such an issue was tackled in \cite{mr13,jkm} for the case of $A_2$ Toda system. The results in these papers exploit the structure of $A_2$ matrix and cannot be extended to other systems. Since in the systems $B_2,G_2$ the interaction is stronger, and asymmetric, we expect a similar approach to be even harder.\\

The same argument also allows to get a multiplicity result when the energy is a Morse functional, which is a generic condition in a sense which will be clarified in Section $2$.\\
In fact, as was done in \cite{bdm,bat1}, if $\Si$'s genus equals $g>1$ we can take as $\g_1,\g_2$ not just two circles but two wedge sums of $g$ circles each. This yields higher-dimensional homology groups for low sub-levels which in turn gives, through Morse theory, a higher number of solutions.\\

The first result we prove is precisely the following:\\

\bthm\label{exmolt}${}$\\
Suppose $\chi(\Si)\le0$ and $\rho_1,\rho_2\nin4\pi\N$. Then, the $B_2$ Toda system \eqref{b2toda} admits solutions.\\
Moreover, if $\rho_1\in(4K_1\pi,4(K_1+1)\pi)\x(4K_2\pi,4(K_2+1)\pi)$, then for a generic choice of $g,h_1,h_2$ the problem has at least $\bin{K_1+\l\lfloor\fr{-\chi(\Si)}2\r\rfloor}{\l\lfloor\fr{-\chi(\Si)}2\r\rfloor}\bin{K_2+\l\lfloor\fr{-\chi(\Si)}2\r\rfloor}{\l\lfloor\fr{-\chi(\Si)}2\r\rfloor}$ solutions.\\
The same results hold true for the $G_2$ Toda system \eqref{g2toda}, provided $\rho_1<4\pi\l(2+\sqrt2\r),\,\rho_2<4\pi\l(5+\sqrt7\r)$.
\ethm\

In the final part of this paper we extend some of the analysis of the sub-levels to the case of general \emph{competitive} Liouville systems with positive singularities.\\
We will consider systems of the type:
\bequ\label{general}
-\D u_i=\sum_{i=1}^Na_{ij}\rho_j\l(\fr{h_je^{u_j}}{\int_\Si h_ie^{u_i}\mrm dV_g}-1\r)-\sum_{m=1}^M\a_{im}(\d_{p_m}-1)\q\q\q i=1,\ds,N,
\eequ
where $A=(a_{ij})_{i,j=1,\ds,N}$ is a symmetric, positive definite matrix with non-positive entries outside the diagonal (i.e., $a_{ij}\le0$ for any $i\ne j$); $p_1,\ds,p_M$ are given points of $\Si$ and $\a_{im}$ are non-negative numbers for $i=1,\ds,N,m=1,\ds,M$.\\
Such systems can be written in a variational form with a simple manipulation involving the Green's function $G_p$ of $-\D$, namely the only solution of
$$\l\{\bll-\D G_p=\d_p-1\\\int_\Si G_p=0\earr\r.:$$
the change of variable
$$u_i\to u_i+\sum_{m=1}^M\a_{im}G_{p_m}$$
transforms problem \eqref{general} in
\bequ\label{system}
-\D u_i=\sum_{i=1}^Na_{ij}\rho_j\l(\fr{\wt h_je^{u_j}}{\int_\Si\wt h_ie^{u_i}\mrm dV_g}-1\r)\q\q\q i=1,\ds,N,
\eequ
where $\wt h_i:=h_ie^{-4\pi\sum_{m=1}^M\a_{im}G_{p_m}}$ has the same behavior as $d(\cd,p_m)^{2\a_{im}}$ around each $p_m$.\\
The energy associated to \eqref{system} is given by
\bequ\label{ja}
J_{A,\rho}(u)=\int_\Si Q_A(u)\mrm dV_g-\sum_{i=1}^N\rho_i\l(\log\int_\Si\wt h_ie^{u_i}\mrm dV_g-\int_\Si u_i\mrm dV_g\r),
\eequ
with $Q_A(u)=\fr{1}2\sum_{i,j=1}^Na^{ij}\n u_i\cd\n u_j$ and $a^{ij}$ are the entries of the inverse matrix $A^{-1}$ of $A$.\\
The functional $J_{A,\rho}$ has been proved in \cite{bat2} to be coercive if and only if $\rho_i<\fr{8\pi}{a_{ii}}$ for all $i=1,\ds,N$.\\
Actually, the requirement for $A$ to be symmetric can be slightly weakened: we may assume the matrix to be \emph{symmetrizable}, namely we can write it as a product $A=SD$ of a symmetric matrix $S$ times a diagonal matrix $D$. Such matrices can be led back to symmetric matrices in the same way as was done before with $B_2$ and $G_2$: roughly speaking, we multiply each parameter $\rho_i$ by the $i^\mrm{th}$ element on the diagonal of $D$. Therefore, we will omit this fact and consider only the case of a symmetric $A$.\\
As in the case of Theorem \ref{exmolt}, we will show that very low sub-levels of the energy functional are not contractible.\\
As before, we will exploit the existence of retractions from $\Si$ to one or more circles. This will prevent some issues related to the interaction between different components $u_i,u_j$ with $i\ne j$, hence the argument will work regardless of the structure of the matrix $A$.\\

\bthm\label{sublev}${}$\\
Suppose $\chi(\Si)\le0$ and $\rho_i\nin\fr{8\pi}{a_{ii}}\N$ for all $i$'s.\\
Then, there exists $L\gg0$ such that $\l\{u\in H^1(\Si)^N;\,J_{A,\rho}(u)\le-L\r\}$ is not contractible.
\ethm\

Anyway, unlike in Theorem \ref{exmolt}, such a result does not suffice to yield existence of solutions.\\
This is because, in the general case, the compactness assumptions which are needed to use this approach, and in particular quantization of local blow-up limits, are not known to hold. This seems to be quite a difficult open problem in general, as it has been solved just for few specials $2$-dimensional systems (see \cite{lwz,lz}).\\
Such a compactness result would imply, for given $A$ and $\rho$, existence of solutions.\\

\bcor\label{exist}${}$\\
Suppose $\ol\rho$ satisfies $\rho_i\nin\fr{8\pi}{a_{ii}}\N$ for all $i$'s and it has a neighborhood $\mcal N\sub\R_{>0}^N$ such that the set of solution $\{u_\rho\}_{\rho\in\mcal N}$ satisfying $\int_\Si u_\rho\mrm dV_g=0$ is compact in $H^1(\Si)^N$.\\
Then, \eqref{system} is solvable for $\rho=\ol\rho$.
\ecor\

The content of the paper is the following. In Section $2$ we provide some notation and preliminary results which will be used throughout the whole paper. In Section $3$ we construct a family of test functions whose image is contained in very low sub-levels of the energy functional. Section $4$ is devoted to the proof of the improved Moser-Trudinger inequality. In Section $5$ we put together the previously obtained results to prove Theorem \ref{exmolt}. Finally, in Section $6$ we adapt the arguments to prove Theorem \ref{sublev}.\\

\section{Notations and preliminaries}\

We will provide here some notation and some known useful preliminary results.\\

The metric distance between two points $x,y\in\Si$ will be denoted by $d(x,y)$; similarly, for $\O,\O'\sub\Si$ we will denote:
$$d(x,\O):=\inf\{d(x,y):\,x\in\O\},\q\q\q d(\O,\O'):=\inf\{d(x,y):\,x\in\O,\,y\in\O'\}.$$
To denote the diameter of a set $\O$ we will write:
$$\mrm{diam}(\O):=\sup\{d(x,y):\,x,y\in\O\}.$$
The open metric ball centered in $p$ having radius $r$ as
$$B_r(x):=\{y\in\Si:\,d(x,y)<r\}.$$
Similarly, for $\O\sub\Si$ we will write
$$B_r(\O):=\{y\in\Si:\,d(y,\O)<r\}.$$\

The standard functional spaces will be denoted as $L^p(\O),C^\infty(\Si),C^\infty(\Si)^N,\ds$. A subscript will be added to denote vectors with positive component or (almost everywhere) positive functions, like $\R_{>0},C^\infty_{>0}(\Si)$.\\
The positive and negative part of a real number $t$ will be denoted respectively by $t^+:=\max\{0,t\}$ and $t^-:=\max\{0,-t\}$.\\

For any continuous map $f:\Si\to\Si$ and any measure $\mu$ defined on $\Si$, we define the push-forward of $\mu$ with respect to $f$ as the measure defined by
$$f_*\mu(B)=\mu\l(f^{-1}(B)\r).$$
We stress that the push-forward of finitely-supported measures has a particularly simple form:
$$\mu=\sum_{k=1}^Kt_k\d_{x_k}\q\q\q\To\q\q\q f_*\mu=\sum_{k=1}^Kt_k\d_{f(x_k)}.$$\

Given a function $u\in L^1(\Si)$ and a measurable set $\O\sub\Si$ with positive measure, we will denote the average of $u$ on $\O$ as
$$\fint_\O u\mrm dV_g=\fr{1}{|\O|}\int_\O u\mrm dV_g.$$
In particular, since $|\Si|=1$, we can write
$$\int_\Si u\mrm dV_g=\fint_\Si u\mrm dV_g.$$
We will indicate the subset of $H^1(\Si)$ containing functions with zero average as
$$\ol H^1(\Si):=\l\{u\in H^1(\Si):\,\int_\Si u=0\r\}.$$
Since the functionals $J_{B_2,\rho},J_{G_2,\rho}$ defined by \eqref{jb2} and \eqref{jg2} are both invariant by addition of constants, as well as the systems \eqref{b2toda} and \eqref{g2toda}, it will not be restrictive to study both of them on $\ol H^1(\Si)^2$ rather than on $H^1(\Si)^2$.\\
As anticipated, the sub-levels of the energy functionals will play an essential role throughout most of the paper. They are denoted as
\bequ\label{subl}
J_{B_2,\rho}^a=\l\{u\in H^1(\Si)^2:\,J_{B_2,\rho}(u)\le a\r\},\q\q\q J_{G_2,\rho}^a=\l\{u\in H^1(\Si)^2:\,J_{G_2,\rho}(u)\le a\r\}.
\eequ\

The composition of two homotopy equivalences $F_1:X\x[0,1]\to Y$ and $F_2:Y\x[0,1]\to Z$ satisfying $F_1(\cd,1)=F_2(\cd,0)$ is the map $F_2\ast F_1:X\x[0,1]\to Z$ defined by
$$F_2\ast F_1:(x,s)\mapsto\l\{\bll F_1(x,2s)&\tx{if }s\le\fr{1}2\\F_2(x,2s-1)&\tx{if }s>\fr{1}2\earr\r..$$
We will denote the identity map on $X$ as $\mrm{Id}_X$.\\
We will denote the $q^\mrm{th}$ homology group with coefficient in $\Z$ of a topological space $X$ as $H_q(X)$. Isomorphisms between homology groups will be denoted just by an equality sign.\\
Reduced homology groups will be denoted as $\wt H_q(X)$, namely
$$H_0(X)=\wt H_0(X)\pl\Z,\q\q\q H_q(X)=\wt H_q(X)\q\tx{if }q\ge1.$$
The $q^\mrm{th}$ Betti number of $X$, that is the dimension of its $q^\mrm{th}$ group of homology, will be indicated by $b_q(X):=\mrm{rank}(H_q(X))$.\\
The symbol $\wt b_q(X)$ will stand for the $q^\mrm{th}$ reduced Betti number, namely the dimension of $\wt H_q(X)$, that is
$$\wt b_0(X)=b_0(X)-1,\q\q\q\wt b_q(X)=b_q(X)\q\tx{if }q\ge1.$$\\
If $J_{B_2,\rho}$ is a Morse function, the symbol $\mcal C_q(J_{B_2};a,b)$ will indicate the number of critical points $u$ of $J_{B_2,\rho}$ with Morse index $q$ satisfying $a\le J_{B_2,\rho}(u)\le b$. The total number of critical points of index $q$ will be denoted as $\mcal C_q(J_{B_2,\rho})$, namely $\mcal C_q(J_{B_2,\rho}):=\mcal C_q(J_{B_2,\rho},+\infty,-\infty)$. A similar notation will be used for $J_{G_2,\rho}$.\\

We will indicate with the letter $C$ large constants, which can vary among different lines and formulas. To underline the dependence of $C$ on some parameter $\a$, we will write $C_\a$ and so on.\\
We will denote as $o_\a(1)$ quantities which tend to $0$ as $\a$ tends to $0$ or to $+\infty$ and we will similarly indicate bounded quantities as $O_\a(1)$, omitting in both cases the subscript(s) when evident from the context.\\

Let us now report the compactness results for solutions of \eqref{b2toda} and \eqref{g2toda}.\\
We start by a concentration-compactness argument from \cite{ln,bat2,batman}.\\

\bthm(\cite{ln}, Theorem $4.2$; \cite{bat2}, Theorem $3.1$; \cite{batman}, Theorem $2.1$)\label{conccomp}${}$\\
Let $\l\{u^n=\l(u_1^n,u_2^n\r)\r\}_{n\in\N}$ be a sequence of solutions of \eqref{b2toda} or \eqref{g2toda} with $\rho_i^n\us{n\to+\infty}\lto\rho_i$. Define
$$S_i:=\l\{x\in\Si:\,\ex\,\l\{x^n\r\}_{n\in\N}\sub\Si:\,u_i^n\l(x^n\r)-\log\int_\Si h_ie^{u_i^n}\mrm dV_g\us{n\to+\infty}\lto+\infty\r\}.$$
Then, up to subsequences, one of the following alternative occurs:
\begin{itemize}
\item(Compactness) For each $i=1,2$ $u_i^n-\log\int_\Si h_ie^{u_i^n}\mrm dV_g$ is uniformly bounded in $L^\infty(\Si)$.
\item(Blow-up) The blow-up set $\mcal S:=\mcal S_1\cup\mcal S_2$ is non-empty and finite.
\end{itemize}
Moreover,
$$\rho_i^n\fr{h_ie^{u_i^n}}{\int_\Si h_ie^{u_i^n}\mrm dV_g}\us{n\to+\infty}\wk r_i+\sum_{x\in\mcal S_i}\s_i(x)\d_x$$
in the sense of measures, with $r_i\in L^1(\O)$ and $\s_i(x)$ defined by
$$\s_i(x):=\lim_{r\to0}\lim_{n\to+\infty}\rho_i^n\fr{\int_{B_r(x)}h_ie^{u_i^n}\mrm dV_g}{\int_\Si h_ie^{u_i^n}\mrm dV_g}.$$
Finally, if $x\in\mcal S_1$ and $2\s_1(x)-\s_2(x)\ge4\pi$, then $r_1\eq0$, whereas if $x\in\mcal S_2$ and $2\s_2(x)+a_{21}\s_2(x)\ge4\pi$, then $r_2\eq0$, with $a_{21}=-2$ in the case of \eqref{b2toda} and $a_{21}=-3$ in the case of \eqref{g2toda}.
\ethm\

We next have a quantization result for the values $\s_1(x),\s_2(x)$:\\

\bthm(\cite{lz})\label{quant}${}$\\
Let $\mcal S,\s_i(x)$ be defined as in Theorem \ref{conccomp} and suppose $x\in\mcal S$.\\
In the case of \eqref{b2toda}, $(\s_1(x),\s_2(x))$ is one of the following:
$$(4\pi,0)\q\q\q(0,4\pi)\q\q\q(4\pi,12\pi)\q\q\q(8\pi,4\pi)\q\q\q(12\pi,12\pi)\q\q\q(8\pi,16\pi)\q\q\q(12\pi,16\pi).$$
In the case of \eqref{g2toda}, under the extra assumption $\rho_1<4\pi\l(2+\sqrt2\r),\,\rho_2<4\pi\l(5+\sqrt7\r)$, $(\s_1(x),\s_2(x))$ is one of the following:
$$(4\pi,0)\q\q\q(0,4\pi)\q\q\q(4\pi,16\pi)\q\q\q(8\pi,4\pi)\q\q\q(8\pi,24\pi)$$
In particular, either $r_1\eq0$ or $r_2\eq0$.
\ethm\

A straightforward applications of Theorems \ref{conccomp} and \ref{quant} gives the following crucial Corollary:\\

\bcor\label{comp}${}$\\
Define $\G:=4\pi\N\x\R\cup\R\x4\pi\N$.\\
The family of solutions $\{u_\rho\}_{\rho\in\mcal K}\sub\ol H^1(\Si)^2$ of \eqref{b2toda} is uniformly bounded in $W^{2,p}(\Si)^2$ for some $p>1$ for any given $\mcal K\Sub\R_{>0}^2\sm\G$.\\
The same holds true for solutions of \eqref{g2toda}, provided $\mcal K\Sub\l(0,4\pi\l(2+\sqrt2\r)\r)\x\l(0,4\pi\l(5+\sqrt7\r)\r)\sm\G$.
\ecor\

The same argument from \cite{luc}, with minor modifications, shows the existence of bounded Palais-Smale sequences. This fact and the compactness result allow to get the following deformation Lemma.\\

\blem\label{deform}${}$\\
Let $a,b\in\R$ be such that $a<b$ and $J_{B_2,\rho}$ (respectively, $J_{G_2,\rho}$) has no critical point with $a\le J_{B_2,\rho}\le b$ (respectively, $a\le J_{B_2,\rho}\le b$).\\
If $\rho\nin\G$, then, $J_{B_2,\rho}^a$ (respectively, $J_{G_2,\rho}^a$) is a deformation retract of $J_{B_2,\rho}^b$ (respectively, of $J_{G_2,\rho}^b$).
\elem\

Compactness of solutions also implies boundedness from above of the energy on solutions, hence the following:\\

\bcor\label{contr}${}$\\
If $\rho\nin\G$, then there exists $L>0$ such that $J_{B_2,\rho}^L$ and $J_{G_2,\rho}^L$ are deformation retracts of $H^1(\Si)^2$. In particular, they are contractible.
\ecor\

Morse inequalities provide an estimate on the number of solutions.\\

\blem\label{morsein}${}$\\
If $\rho\nin\G$ and $J_{B_2,\rho}$ and $J_{G_2,\rho}$ are Morse functions, then there exist $L\gg0$ such that $\mcal C_q(J_{B_2,\rho};-L,L)\ge\wt b_q\l(J_{B_2,\rho}^{-L}\r)$ and $\mcal C_q(J_{G_2,\rho};-L,L)\ge\wt b_q\l(J_{G_2,\rho}^{-L}\r)$.\\
In particular,
$$\#\tx{ solutions of }\eqref{b2toda}=\sum_{q=0}^{+\infty}\mcal C_q\l(J_{B_2,\rho}\r)\ge\sum_{q=0}^{+\infty}\mcal C_q\l(J_{B_2,\rho};-L,L\r)\ge\sum_{q=0}^{+\infty}\l(\wt b_q\l(J_{B_2,\rho}^{-L}\r)-\wt b_q\l(J_{B_2,\rho}^L\r)\r)=\sum_{q=0}^{+\infty}\wt b_q\l(J_{B_2,\rho}^{-L}\r),$$
$$\#\tx{ solutions of }\eqref{g2toda}=\sum_{q=0}^{+\infty}\mcal C_q\l(J_{G_2,\rho}\r)\ge\sum_{q=0}^{+\infty}\mcal C_q\l(J_{G_2,\rho};-L,L\r)\ge\sum_{q=0}^{+\infty}\l(\wt b_q\l(J_{G_2,\rho}^{-L}\r)-\wt b_q\l(J_{G_2,\rho}^L\r)\r)=\sum_{q=0}^{+\infty}\wt b_q\l(J_{G_2,\rho}^{-L}\r).$$
\elem\

By arguing as in \cite{bdm,bat1}, we find that $J_{B_2,\rho}$ and $J_{G_2,\rho}$ are Morse functionals for a generic choice of initial data:\\

\blem\label{dense}${}$\\
There exists a dense open set $\mcal D\sub\mcal M^2(\Si)\x C(\Si)^2$ such that for any $(g,h_1,h_2)\in\mcal D$ both $J_{B_2,\rho}$ and $J_{G_2,\rho}$ are Morse functions.
\elem\

Here are the Moser-Trudinger inequalities we will need.\\
The first is a very classic result from \cite{mos,fon}:\\

\bthm(\cite{mos}, Theorem $2$; \cite{fon}, Theorem $1.7$)\label{mtscal}\\
There exists $C>0$ such that for any $u\in H^1(\Si)$ one has
$$16\pi\l(\log\int_\Si e^u\mrm dV_g-\int_\Si u\mrm dV_g\r)\le\int_\Si|\n u|^2\mrm dV_g+C.$$
\ethm\

We also need a vectorial generalization of Theorem \ref{mtscal}, which was given in \cite{bat2} in a pretty much general form.\\
Notice that, although its original form concerns symmetric matrices, it actually works even for the $B_2$ and $G_2$ matrices just by writing
$$\l(\barr{cc}2&-1\\-2&2\earr\r)\l(\barr{c}\rho_1\\\rho_2\earr\r)=\l(\barr{cc}2&-2\\-2&4\earr\r)\l(\barr{c}\rho_1\\\fr{\rho_2}2\earr\r)\q\q\q\l(\barr{cc}2&-1\\-3&2\earr\r)\l(\barr{c}\rho_1\\\rho_2\earr\r)=\l(\barr{cc}2&-3\\-3&6\earr\r)\l(\barr{c}\rho_1\\\fr{\rho_2}3\earr\r).$$\

\bthm(\cite{bat2}, Theorem $1.2$)\label{mt}\\
There exists $C>0$ such that for any $u\in H^1(\Si)^2$ one has
$$4\pi\l(\log\int_\Si e^{u_1}\mrm dV_g-\int_\Si u_1\mrm dV_g\r)+2\pi\l(\log\int_\Si e^{u_2}\mrm dV_g-\int_\Si u_2\mrm dV_g\r)\le\int_\Si Q_{B_2}(u)\mrm dV_g+C$$
$$4\pi\l(\log\int_\Si e^{u_1}\mrm dV_g-\int_\Si u_1\mrm dV_g\r)+\fr{4}3\pi\l(\log\int_\Si e^{u_2}\mrm dV_g-\int_\Si u_2\mrm dV_g\r)\le\int_\Si Q_{G_2}(u)\mrm dV_g+C$$
In particular, both \eqref{jb2} and \eqref{jg2} are bounded from below if and only if $\rho_1,\rho_2\le4\pi$.\\
Moreover, they are coercive if and only if $\rho_1,\rho_2<4\pi$. In this case, both \eqref{b2toda} and \eqref{g2toda} have a minimizing solution.
\ethm\

The following is a covering lemma which we will need to prove the improved Moser-Trudinger inequality.\\

\blem(\cite{bjmr}, Lemma $4.1$)\label{cover}\\
Let $\d>0,K_1,K_2\in\N$ be given, $f_1,f_2\in L^1(\Si)$ be positive a.e. and such that $\|f_1\|_{L^1(\Si)}=\|f_2\|_{L^1(\Si)}=1$ and $\{\O_{1k}\}_{k=1}^{K_1},\{\O_{2k}\}_{k=1}^{K_2}$ satisfy
\beqy
d(\O_{ik},\O_{ik'})\ge\d&\q\q\q&\fa\,i=1,2,\,k,k'=1,\ds,K_i,\,k\ne k',\\
\int_{\O_{ik}}f_i\mrm dV_g\ge\d&\q\q\q&\fa\,i=1,2,\,k=1,\ds,K_i.
\eeqy
Then, there exists $\d'=\d'(\d,\Si)$ and $\{\O_k\}_{k=1}^{\max\{K_1,K_2\}}$ such that
\beqy
d(\O_k,\O_{k'})\ge\d'&\q\q\q&\fa\,k,k'=1,\ds,\max\{K_1,K_2\},\,k\ne k',\\
\int_{\O_k}f_i\mrm dV_g\ge\d'&\q\q\q&\fa\,i=1,2,\,k=1,\ds,K_i.
\eeqy
\elem\

As anticipated, the assumption on $\Si$ yields a simple but very powerful result from general topology:\\

\blem\label{retra}${}$\\
Let $\Si$ be a surface of genus $g=\l\lfloor\fr{-\chi(\Si)}2\r\rfloor\ge1$.\\
Then, there exist two disjointed curves $\g_1,\g_2$, each of which is homeomorphic to a wedge sum of $g$ circles, and two global retractions $\Pi_i:\Si\to\g_i$.
\elem\

Let us now introduce the space of \emph{barycenters} on a measure space $X$, that is the space of measures supported in at most $K$-points:
\bequ\label{xk}
(X)_K:=\l\{\sum_{k=1}^Kt_k\d_{x_k}:\,x_k\in X,\,t_k\ge0,\,\sum_{k=1}^Kt_k=1\r\}
\eequ
Such a space will be endowed with tha $Lip'$ norm, namely the norm defined by duality with the space $\mrm{Lip}(\Si)$:
$$\|\mu\|_{\mrm{Lip}'(\Si)}:=\sup_{\phi\in\mrm{Lip}(\Si),\|\phi\|_{\mrm{Lip}(\Si)}\le1}\l|\int_\Si\phi\mrm d\mu\r|.$$

A first important property of such a space is being a Euclidean neighborhood retract:\\

\blem(\cite{dm}, Lemma $3.10$)\label{psik}\\
Let $(\Si)_K$ be defined by \eqref{xk}.\\
Then, there exist $\e_0>0$ and a retraction
$$\psi_K:\{\mu\in\mcal M(\Si):\,d(\mu,(\Si)_K)<\e_0\}\to(\Si)_K.$$
In particular, if $\mu^n\us{n\to+\infty}\wk\s$ for some $\s\in(\Si)_K$, then $\psi_K(\mu^n)\us{n\to+\infty}\lto\s$.
\elem\

We then introduce the \emph{join} of two topological spaces $X$ and $Y$.\\
It is basically the product of the two spaces and the unit interval, with two identifications made at the endpoints: when $t=0$ the space $Y$ is collapsed and when $t=1$ the space $X$ is collapsed.\\
Precisely, we define:
$$X\star Y:=\fr{X\x Y\x[0,1]}\sim,$$
where $\sim$ is given by
$$(x,y,0)\sim(x,y',0)\q\fa\,x\in X,\,\fa\,y,y'\in Y\q\q\q\q\q\q(x,y,1)\sim(x',y,1)\q\fa\,x,x'\in X,\,\fa\,y\in Y.$$
Such an object is used as a model for a (non-exclusive) alternative between $X$ and $Y$.\\
In this paper, we will use the join to express the following rough idea: as $J_{B_2,\rho}$ (respectively, $J_{G_2,\rho}$) is lower and lower, at least one of the unit measures
\bequ\label{fiu}
f_{i,u}=\fr{h_ie^{u_i}}{\int_\Si h_ie^{u_i}\mrm dV_g}
\eequ
is almost supported at a finite number of points. To express the concentration at at most $K$ points we will use the space of barycenters on $\Si$ or on $\g_i$, whereas the join will express the alternative between the two components.\\
The homology of such a join is expressed by this proposition:\\

\bprop(\cite{bdm}, Proposition $3.2$; \cite{hat}, Theorem $3.21$)\label{hom}\\
Let $\g_1,\g_2$ be wedge sums of $g$ circles.\\
Then,
$$H_q((\g_1)_{K_1}\star(\g_2)_{K_2})=\l\{\bll\Z&\tx{if }q=0\\\Z^{\bin{K_1+g-1}{g-1}\bin{K_2+g-1}{g-1}}&\tx{if }q=2K_1+2K_2-1\\0&\tx{if }q\ne0,2K_1+2K_2-1\earr\r.$$ 
\eprop\

\section{Test functions}\

In this section we will consider two families of test functions, modeled on $(\g_1)_{K_1}\star(\g_2)_{K_2}$, on which $J_{B_2,\rho}$ and $J_{G_2,\rho}$ respectively attain arbitrarily low values.\\
In other words, we will get a family of maps from the join of the barycenters' space to arbitrarily low sub-levels.

\bthm\label{testfun}${}$\\
Define, for any $\la>0$ and $\z=(\s_1,\s_2,t)=\l(\sum_{k=1}^{K_1}t_{1k}\d_{x_{1k}},\sum_{k=1}^{K_2}t_{2k}\d_{x_{2k}},t\r)\in(\g_1)_{K_1}\star(\g_2)_{K_2}$,
$$\ph_1=\ph_1^\la(\z)=\log\sum_{k=1}^{K_1}\fr{t_{1k}}{\l(1+(\la(1-t)d(\cd,x_{1k}))^2\r)^2}\q\q\q\ph_2=\ph_2^\la(\z)=\log\sum_{k=1}^{K_2}\fr{t_{2k}}{\l(1+(\la td(\cd,x_{2k}))^2\r)^2}$$
$$\Phi_{B_2}^\la(\z)=\l(\ph_1-\fr{\ph_2}2,\ph_2-\ph_1\r)\q\q\q\Phi_{G_2}^\la(\z)=\l(\ph_1-\fr{\ph_2}2,\ph_2-\fr{3}2\ph_1\r).$$
If $\rho\in(4K_1\pi,4(K_1+1)\pi)\x(4K_2\pi,4(K_2+1)\pi)$, then
$$\l.\bl J_{B_2,\rho}\l(\Phi_{B_2}^\la(\z)\r)\us{\la\to+\infty}\lto-\infty\\J_{G_2,\rho}\l(\Phi_{G_2}^\la(\z)\r)\us{\la\to+\infty}\lto-\infty\earr\r.\q\q\q\tx{uniformly in }\z\in(\g_1)_{K_1}\star(\g_2)_{K_2}.$$
\ethm\

The proof of Theorem \ref{testfun} will follow immediately by the following three lemmas, where the different parts of the functional $J_{B_2,\rho},J_{G_2,\rho}$ are estimated separately:\\

\blem\label{grad}${}$\\
Let $\z,\ph_i$ be as in Theorem \ref{testfun}. Then,
$$\int_\Si Q_{B_2}\l(\ph_1-\fr{\ph_2}2,\ph_2-\ph_1\r)\mrm dV_g\le8K_1\pi\log\max\{1,\la(1-t)\}+4K_2\pi\log\max\{1,\la t\}+C$$
$$\int_\Si Q_{G_2}\l(\ph_1-\fr{\ph_2}2,\ph_2-\fr{3}2\ph_1\r)\mrm dV_g\le8K_1\pi\log\max\{1,\la(1-t)\}+\fr{8}3K_2\pi\log\max\{1,\la t\}+C$$
\elem\

\bpf${}$\\

Since we can write

$$\int_\Si Q_{B_2}\l(\ph_1-\fr{\ph_2}2,\ph_2-\ph_1\r)\mrm dV_g=\fr{1}4\int_\Si|\n\ph_1|^2\mrm dV_g-\fr{1}4\int_\Si\n\ph_1\cd\n\ph_2\mrm dV_g+\fr{1}8\int_\Si|\n\ph_2|^2\mrm dV_g$$
$$\int_\Si Q_{G_2}\l(\ph_1-\fr{\ph_2}2,\ph_2-\fr{3}2\ph_1\r)\mrm dV_g=\fr{1}4\int_\Si|\n\ph_1|^2\mrm dV_g-\fr{1}4\int_\Si\n\ph_1\cd\n\ph_2\mrm dV_g+\fr{1}{12}\int_\Si|\n\ph_2|^2\mrm dV_g,$$
then we will suffice to show
\bequ\label{grad1}
\int_\Si|\n\ph_1|^2\mrm dV_g\le32K_1\pi\log\max\{1,\la(1-t)\}+C
\eequ
\bequ\label{grad12}
\int_\Si\n\ph_1\cd\n\ph_2\mrm dV_g=O(1)
\eequ
\bequ\label{grad2}
\int_\Si|\n\ph_2|^2\mrm dV_g\le32K_2\pi\log\max\{1,\la t\}+C.
\eequ
Since $|\n d(\cd,x_{1k})|=1$ a.e. on $\Si$, it holds
\beqy
|\n\ph_1|&=&\l|\fr{\sum_{k=1}^{K_1}\fr{-4t_{1k}(\la(1-t))^2d(\cd,x_{1k})\n d(\cd,x_{1k})}{\l(1+\l(\la(1-t)d\l(\cd,x_{1k}\r)\r)^2\r)^2}}{\sum_{k=1}^{K_1}\fr{t_{1k}}{1+(\la(1-t)d(\cd,x_{1k}))^2}}\r|\\
&\le&\fr{\sum_{k=1}^{K_1}\fr{4t_{1k}(\la(1-t))^2d(\cd,x_{1k})}{\l(1+\l(\la(1-t)d\l(\cd,x_{1k}\r)\r)^2\r)^2}}{\sum_{k=1}^{K_1}\fr{t_{1k}}{1+(\la(1-t)d(\cd,x_{1k}))^3}}\\
&\le&\max_k\fr{4(\la(1-t))^2d(\cd,x_{1k})}{\l(1+(\la(1-t)d(\cd,x_{1k}))^2\r)^2}\\
&\le&\min\l\{2\la(1-t),\fr{4}{\min_kd(\cd,x_{1k})}\r\}.
\eeqy
In the same way, we get the estimate $|\n\ph_2|\le\min\l\{2\la t,\fr{4}{\min_kd(\cd,x_{2k})}\r\}$.\\

In view of such estimates, we divide $\Si$ in $K_1$ regions $\O_{11},\ds,\O_{1K_1}$, depending on which of the points $x_{ik}$ is the closest:
$$\O_{1k}:=\l\{x\in\Si:\,d(x,x_{1k})=\min_{k'}d(x,x_{1k'})\r\};$$
we similarly divide $\Si$ in regions $\O_{21},\ds,\O_{2K_2}$.\\
By splitting the integral in such regions, we get
\beqy
\int_\Si|\n\ph_1|^2\mrm dV_g&\le&\sum_{k=1}^{K_1}\int_{\O_{1k}}\min\l\{4(\la(1-t))^2,\fr{16}{d(\cd,x_{1k})^2}\r\}\\
&\le&\sum_{k=1}^{K_1}\l(4(\la(1-t))^2\int_{B_\fr{1}{\la(1-t)}(x_{1k})}\mrm dV_g+16\int_{\Si\sm B_\fr{1}{\la(1-t)}(x_{1k})}\fr{\mrm dV_g}{d(\cd,x_{1k})}\r)\\
&\le&\sum_{k=1}^{K_1}(C+32\pi\log\max\{1,\la(1-t)\})\\
&\le&32K_1\pi\log\max\{1,\la(1-t)\}+C.
\eeqy
This proves \eqref{grad1}, whereas \eqref{grad2} can be proved in the same way.\\

To obtain \eqref{grad12} we exploit the distance between the points $x_{1k}$ and $x_{2k'}$, which is bounded from below by a positive constant.\\
Taking $\d=\fr{d(\g_1,\g_2)}2$ one gets $B_\d(x_{1k})\cap B_\d(x_{2k'})=\es$ for any $k,k'$, therefore:
\beqy
\l|\int_\Si\n\ph_1\cd\n\ph_2\mrm dV_g\r|&\le&\sum_{k,k'}\int_{\O_{1k}\cap\O_{2k'}}|\n\ph_1||\n\ph_2|\mrm dV_g\\
&\le&\sum_{k,k'}\int_{\O_{1k}\cap\O_{2k'}}\fr{\mrm dV_g}{d(\cd,x_{1k})d(\cd,x_{2k'})}\\
&\le&4\sum_{k,k'}\l(\int_{\O_{1k}\cap\O_{2k'}\sm B_\d(x_{1k})}\fr{\mrm dV_g}{\d d(\cd,x_{2k'})}+\int_{\O_{1k}\cap\O_{2k'}\sm B_\d(x_{2k'})}\fr{\mrm dV_g}{\d d(\cd,x_{1k})}\r)\\
&\le&\fr{4}\d\sum_{k,k'}\int_\Si\l(\fr{1}{d(\cd,x_{2k'})}+\fr{1}{d(\cd,x_{1k})}\r)\mrm dV_g\\
&\le&C.
\eeqy

\epf\

\blem\label{average}${}$\\
Let $\z,\ph_i$ be as in Theorem \ref{testfun}. Then,
$$\int_\Si\ph_1\mrm dV_g=-4\log\max\{1,\la(1-t)\}+O(1),\q\q\q\int_\Si\ph_2\mrm dV_g=-4\log\max\{1,\la t\}+O(1).$$
\elem\

\bpf${}$\\
The proof will be given only for $i=1$, since the other estimate is similar.\\
We first notice that
$$\ph_i=4\log\sum_{k=1}^{K_1}\fr{t_{1k}}{\max\{1,\la(1-t)d(\cd,x_{1k})\}}+O(1).$$
Then, we can provide an estimate from above:
\beqy
\int_\Si\ph_1\mrm dV_g&\le&4\int_\Si\log\min\l\{0,\fr{1}{\la(1-t)\min_kd(\cd,x_{1k})}\r\}\mrm dV_g+O(1)\\
&=&4\int_{\Si\sm\Cup_{k'=1}^{K_1}B_\fr{1}{\la(1-t)}(x_{1k'})}\log\fr{1}{\min_kd(\cd,x_{1k})}\mrm dV_g\\
&-&4\log\max\{1,\la(1-t)\}\l|\Si\sm\Cup_{k'=1}^{K_1}B_\fr{1}{\la(1-t)}(x_{1k'})\r|+O(1)\\
&\le&-4\log\max\{1,\la(1-t)\}+O(1).
\eeqy
The estimate from below is similar:
\beqy
\int_\Si\ph_1\mrm dV_g&\ge&4\int_\Si\log\min\l\{0,\fr{1}{\la(1-t)\max_kd(\cd,x_{1k})}\r\}\mrm dV_g+O(1)\\
&=&4\int_{\Si\sm\Cap_{k'=1}^{K_1}B_\fr{1}{\la(1-t)}(x_{1k'})}\log\fr{1}{\max_kd(\cd,x_{1k})}\mrm dV_g\\
&-&4\log\max\{1,\la(1-t)\}\l|\Si\sm\Cap_{k'=1}^{K_1}B_\fr{1}{\la(1-t)}(x_{1k'})\r|+O(1)\\
&\ge&-4\log\max\{1,\la(1-t)\}+O(1).
\eeqy

\epf\

\blem\label{exp}${}$\\
Let $\z,\ph_i$ be as in Theorem \ref{testfun}. Then,
$$\log\int_\Si h_1e^{\ph_1-\fr{\ph_2}2}\mrm dV_g=-2\log\max\{1,\la(1-t)\}+2\log\max\{1,\la t\}+O(1).$$
$$\log\int_\Si h_2e^{\ph_2-\ph_1}\mrm dV_g=-2\log\max\{1,\la t\}+4\log\max\{1,\la(1-t)\}+O(1),$$
$$\log\int_\Si h_2e^{\ph_2-\fr{3}2\ph_1}\mrm dV_g=-2\log\max\{1,\la t\}+6\log\max\{1,\la(1-t)\}+O(1).$$
\elem\

\bpf${}$\\
As before, we will just prove the first assertion because the same arguments also work for the other two.\\
To give an upper estimates, we write:
\beqy
\int_\Si h_1e^{\ph_1-\fr{\ph_2}2}\mrm dV_g&=&\int_\Si h_1\fr{\sum_{k=1}^{K_1}\fr{t_{1k}}{\l(1+\l(\la(1-t)d\l(\cd,x_{1k}\r)\r)^2\r)^2}}{\sum_{k'=1}^{K_2}\fr{t_{2k'}}{1+\l(\la td\l(\cd,x_{2k'}\r)\r)^2}}\mrm dV_g\\
&\le&C\l(1+(\la t\mrm{diam}(\Si))^2\r)\int_\Si\sum_{k=1}^{K_1}\fr{t_{1k}}{\l(1+\l(\la(1-t)d\l(\cd,x_{1k}\r)\r)^2\r)^2}\mrm dV_g\\
&\le&C\max\{1,\la t\}^2\sum_{k=1}^{K_1}t_{1k}\int_\Si\fr{\mrm dV_g}{\l(1+\l(\la(1-t)d\l(\cd,x_{1k}\r)\r)^2\r)^2}.
\eeqy
Taking, as before, $\d=\fr{d(\g_1,\g_2)}2$, we get a similar estimate from below
\beqy
\int_\Si h_1e^{\ph_1-\fr{\ph_2}2}\mrm dV_g&\ge&\fr{1}C\int_{B_\d(x_{1k})}\fr{\sum_{k=1}^{K_1}\fr{t_{1k}}{\l(1+\l(\la(1-t)d\l(\cd,x_{1k}\r)\r)^2\r)^2}}{\sum_{k'=1}^{K_2}\fr{t_{2k'}}{1+\l(\la td\l(\cd,x_{2k}\r)\r)^2}}\mrm dV_g\\
&\ge&\fr{1}C\l(1+(\la t\d)^2\r)\sum_{k=1}^{K_1}t_{1k}\int_{B_\d(x_{1k})}\fr{\mrm dV_g}{\l(1+\l(\la(1-t)d\l(\cd,x_{1k}\r)\r)^2\r)^2}\\
&\ge&\fr{1}C\max\{1,\la t\}^2\sum_{k=1}^{K_1}t_{1k}\int_{B_\d(x_{1k})}\fr{\mrm dV_g}{\l(1+\l(\la(1-t)d\l(\cd,x_{1k}\r)\r)^2\r)^2}.
\eeqy
Therefore, to conclude the proof, we need to show:
\bequ\label{outball}
\int_{\Si\sm B_\d(x_{1k})}\fr{\mrm dV_g}{\l(1+\l(\la(1-t)d\l(\cd,x_{1k}\r)\r)^2\r)^2}\le C\max\{1,\la t\}^2
\eequ
\bequ\label{inball}
\int_{B_\d(x_{1k})}\fr{\mrm dV_g}{\l(1+\l(\la(1-t)d\l(\cd,x_{1k}\r)\r)^2\r)^2}\sim C\max\{1,\la t\}^2.
\eequ
\eqref{outball} follows from
\beqy
\int_{\Si\sm B_\d(x_{1k})}\fr{\mrm dV_g}{\l(1+\l(\la(1-t)d\l(\cd,x_{1k}\r)\r)^2\r)^2}&\le&\int_{\Si\sm B_\d(x_{1k})}\fr{\mrm dV_g}{\l(1+\l(\la(1-t)\d\r)^2\r)^2}\\
&\le&\fr{C}{\max\{1,\la t\}^4}\\
&\le&\fr{C}{\max\{1,\la t\}^2}.
\eeqy
On the other hand, \eqref{inball} can be obtained through normal coordinates and a change of variables:
\beqy
\int_{B_\d(x_{1k})}\fr{\mrm dV_g}{\l(1+\l(\la(1-t)d\l(\cd,x_{1k}\r)\r)^2\r)^2}&\sim&\int_{B_\d(0)}\fr{\mrm dx}{\l(1+(\la(1-t)|x|)^2\r)^2}\\
&\sim&\fr{1}{(\la(1-t))^2}\int_{B_{\la(1-t)\d}(0)}\fr{\mrm dy}{\l(1+|y|^2\r)^2}\\
&\sim&\fr{1}{\max\{1,\la(1-t)\}^2}.
\eeqy
\epf\

\section{Improved Moser-Trudinger inequality}\

Here, we will prove the following Theorem, which gives important information about low sub-levels.\\

\bthm\label{close}${}$\\
Suppose $\rho\in(4K_1\pi,4(K_1+1)\pi)\x(4K_2\pi,4(K_2+1)\pi)$.\\
Then for any $\e>0$ there exists $L=L_\e>0$ such that any $u\in J_{B_2,\rho}^{-L}\cup J_{G_2,\rho}^{-L}$ verifies, for some $i=1,2$,
$$d_{\mrm{Lip}'}(f_{i,u},(\Si)_{K_i})\le\e,$$
where $f_{i,u}$ is defined by \eqref{fiu}.
\ethm\

The proof of Theorem \ref{close} can be deduced by arguing as in \cite{bjmr,bat1}, Section $4$.\\
We present in details only the main new ingredient, which is a so-called improved Moser-Trudinger inequality.\\
Basically, the constants $4\pi,2\pi,\fr{4}3\pi$ in Theorem \ref{mt} can be multiplied by a suitable integer number under a condition of \qm{spreading} on $f_{1,u},\,f_{2,u}$.\\

Arguing as in \cite{bjmr,bat1}, we can infer by the following Lemma that if $J_{B_2,\rho}(u)\ll0$ (respectively, $J_{G_2,\rho}(u)\ll0$) then either $f_{1,u}$ or $f_{2,u}$ can accumulate mass only around a fixed number of points, hence it must be close to the corresponding space $(\Si)_{K_i}$.\\

\blem\label{improved}${}$\\
Let $\d>0,K_1,K_2\in\N,\,\{\O_{1k}\}_{k=1}^{K_1},\{\O_{2k}\}_{k=1}^{K_2}$ satisfy
\beqy
d(\O_{ik},\O_{ik'})\ge\d&\q\q\q&\fa\,i=1,2,\,k,k'=1,\ds,K_i,\,k\ne k',\\
\int_{\O_{ik}}f_{i,u}\mrm dV_g\ge\d&\q\q\q&\fa\,i=1,2,\,k=1,\ds,K_i.
\eeqy
Then, for any $\e>0$ there exists $C>0$, not depending on $u$, such that
\bequ\label{imprb2}
4\pi K_1\l(\log\int_\Si e^{u_1}\mrm dV_g-\int_\Si u_1\mrm dV_g\r)+2\pi K_2\l(\log\int_\Si e^{u_2}\mrm dV_g-\int_\Si u_2\mrm dV_g\r)\le(1+\e)\int_\Si Q_{B_2}(u)\mrm dV_g+C
\eequ
\bequ\label{imprg2}
4\pi K_1\l(\log\int_\Si e^{u_1}\mrm dV_g-\int_\Si u_1\mrm dV_g\r)+\fr{4}3\pi K_2\l(\log\int_\Si e^{u_2}\mrm dV_g-\int_\Si u_2\mrm dV_g\r)\le(1+\e)\int_\Si Q_{G_2}(u)\mrm dV_g+C
\eequ
\elem\

\bpf${}$\\
Let us assume $K_1\ge K_2$.\\
As a first thing, we apply Lemma \ref{cover} to $f_{1,u},f_{2,u}$ and we get $\{\O_k\}_{k=1}^{K_1}$ such that
$$d(\O_k,\O_{k'})\ge\d'\q\q\q\int_{\O_k}f_{i,u}\mrm dV_g\ge\d'.$$
We then consider, for $k=1,\ds,K_1$, some cut-off functions $\chi_k$ satisfying
\bequ\label{chi}
0\le\chi_k\le1\q\q\q\chi_k|_{\O_k}\eq1\q\q\q\chi_k|_{\Si\sm\O'_k}\eq0\q\tx{with }\O'=B_\fr{\d'}2(\O_k)\q\q\q|\n\chi_k|\le C_{\d'}.
\eequ
Write now $u-\int_\Si u_i\mrm dV_g=v_i+w_i$, with $\int_\Si v_i\mrm dV_g=\int_\Si w_i\mrm dV_g=0$ and $v_i\in L^\infty(\Si)$ (which will be chosen later). It holds:
\beqa
\nonumber\log\int_\Si e^{u_i}\mrm dV_g-\int_\Si u_i\mrm dV_g&\le&\log\int_{\O_k}e^{u_i}\mrm dV_g-\int_\Si u_i\mrm dV_g+\log\fr{1}{\d'}\\
\nonumber&=&\log\int_{\O_k}e^{v_i+w_i}\mrm dV_g+\log\fr{1}{\d'}\\
\nonumber&\le&\log\int_{\O_k}e^{w_i}\mrm dV_g+\|v_i\|_{L^\infty(\Si)}+C\\
\label{uvw}&\le&\log\int_\Si e^{\chi_kw_i}\mrm dV_g+\|v_i\|_{L^\infty(\Si)}+C.
\eeqa
We would like to apply Moser-Trudinger inequality to $\chi_kw_i$. To this purpose, we write
\beqy
\int_\Si|\n(\chi_kw_i)|^2\mrm dV_g&=&\int_\Si\l(\chi_k^2|\n w_i|^2+2(\chi_k\n w_i)\cd(w_i\n\chi_k)+w_i^2|\n\chi_k|^2\r)\mrm dV_g\\
&\le&\int_\Si\l((1+\e)\chi_k^2|\n w_1|^2+\l(1+\fr{1}\e\r)w_i^2|\n\chi_k|^2\r)\mrm dV_g\\
&\le&(1+\e)\int_{\O'_k}|\n w_i|^2\mrm dV_g+C_{\e,\d'}\int_{\O'_k}w_i^2\mrm dV_g.
\eeqy
Using the elementary equalities
\bequ\label{xy}
\fr{|x|^2}2+\fr{x\cd y}2+\fr{|y|^2}4=\fr{|x|^2}4+\fr{|x+y|^2}4\q\q\q|x|^2+x\cd y+\fr{|y|^2}3=\fr{|x|^2}4+\fr{|3x+2y|^2}{12},
\eequ
we similarly get
\bequ\label{qw}
\int_\Si Q_{B_2}(\chi_kw)\mrm dV_g\le(1+\e)\int_{\O'_k}Q_{B_2}(w)\mrm dV_g+C_{\e,\d'}\int_{\O'_k}\l(\fr{w_1^2}2+\fr{w_1w_2}2+\fr{w_2^2}4\r)\mrm dV_g
\eequ
$$\int_\Si Q_{G_2}(\chi_kw)\mrm dV_g\le(1+\e)\int_{\O'_k}Q_{G_2}(w)\mrm dV_g+C_{\e,\d'}\int_{\O'_k}\l(w_1^2+w_1w_2+\fr{w_2^2}3\r)\mrm dV_g$$
Now we choose $v_i$ in order to control the $L^2$ norm of $w_i$. Fixing an orthonormal frame $\{\ph^n\}_{n=1}^{+\infty}$ of eigenfunctions of $-\D$ on $\ol H^1(\Si)$ with positive non-decreasing eigenvalues $\{\la^n\}_{n=1}^{+\infty}$, we write
$$u_i:=\int_\Si u_i\mrm dV_g+\sum_{n=1}^{+\infty}u_i^n\ph^n\q\q\q v_i:=\sum_{n=1}^Nu_i^n\ph^n\q\q\q w^n=\sum_{n=N+1}^{+\infty}u_i^n\ph^n$$
$$\tx{where}\q\q\q N=N_{\e,\d'}:=\max\l\{n\in\N:\la^n<\fr{C_{\e,\d'}}\e\r\}.$$
This choice yields
$$C_{\e,\d'}\int_\Si w_i^2\mrm dV_g\le\e\int_\Si|\n w_i|^2\mrm dV_g\le\e\int_\Si|\n u_i|^2\mrm dV_g$$
and, by \eqref{xy},
$$C_{\e,\d'}\int_\Si\l(\fr{w_1^2}2+\fr{w_1w_2}2+\fr{w_2^2}4\r)\mrm dV_g\le\e\int_\Si Q_{B_2}(w)\mrm dV_g\le\e\int_\Si Q_{B_2}(u)\mrm dV_g$$
$$C_{\e,\d'}\int_\Si\l(w_1^2+w_1w_2+\fr{w_2^2}3\r)\mrm dV_g\le\e\int_\Si Q_{G_2}(w)\mrm dV_g\le\e\int_\Si Q_{G_2}(u)\mrm dV_g.$$
Moreover, we have
\bequ\label{wi}
\int_\Si|w_i|\mrm dV_g\le C\sqrt{\int_\Si|\n w_i|^2\mrm dV_g}\le\e\min\l\{\int_\Si Q_{B_2}(u)\mrm dV_g,\int_\Si Q_{G_2}(u)\mrm dV_g\r\}+C_\e
\eequ
and, since $v_i$ is taken in a finite-dimensional space,
\bequ\label{vi}
\|v_i\|_{L^\infty(\Si)}\le C_N\sqrt{\int_\Si|\n v_i|^2\mrm dV_g}\le\e\min\l\{\int_\Si Q_{B_2}(u)\mrm dV_g,\int_\Si Q_{G_2}(u)\mrm dV_g\r\}+C_\e.
\eequ\

At this point, for $k=1,\ds,K_2$ we apply Theorem \ref{mt} to $\chi_kw$: using the first inequality we get, by \eqref{qw} and \eqref{wi},
\beqy
4\pi\log\int_\Si e^{\chi_kw_1}\mrm dV_g+2\pi\log\int_\Si e^{\chi_kw_2}&\le&4\pi\int_\Si\chi_kw_1\mrm dV_g+2\pi\int_\Si\chi_kw_2\mrm dV_g+\int_\Si Q_{B_2}(\chi_kw)\mrm dV_g+C\\
&\le&4\pi\|\chi_k\|_{L^\infty(\Si)}\int_\Si|w_1|\mrm dV_g+2\pi\|\chi_k\|_{L^\infty(\Si)}\int_\Si|w_2|\mrm dV_g\\
&+&(1+\e)\int_{\O'_k} Q_{B_2}(w)\mrm dV_g+C_{\e,\d'}\int_\Si\l(w_1^2+w_1w_2+\fr{w_2^2}2\r)\mrm dV_g+C\\
&\le&(1+\e)\int_{\O'_k} Q_{B_2}(w)\mrm dV_g+\e'\int_\Si Q_{B_2}(u)\mrm dV_g+C\\
&\le&\int_{\O'_k} Q_{B_2}(w)\mrm dV_g+\e''\int_\Si Q_{B_2}(u)\mrm dV_g+C
\eeqy
Applying now \eqref{uvw} and \eqref{vi} we deduce
\beqa
\nonumber&&4\pi\l(\log\int_\Si e^{u_1}\mrm dV_g-\int_\Si u_1\mrm dV_g\r)+2\pi\l(\log\int_\Si e^{u_2}\mrm dV_g-\int_\Si u_2\mrm dV_g\r)\\
\nonumber&\le&4\pi\log\int_\Si e^{\chi_kw_1}\mrm dV_g+2\pi\log\int_\Si e^{\chi_kw_2}+4\pi\|v_1\|_{L^\infty(\Si)}+2\pi\|v_2\|_{L^\infty(\Si)}+C\\
\label{u1u2}&\le&\int_{\O'_k} Q_{B_2}(w)\mrm dV_g+\e\int_\Si Q_{B_2}(u)\mrm dV_g+C.
\eeqa
For $k=K_2+1,\ds,K_1$ we only have estimates for $f_{1,u}$ on $\O_k$, so we apply the scalar Moser-Trudinger inequality, that is Theorem \ref{mtscal}. By \eqref{xy} we get again the integral of $Q_{B_2}$, hence we can argue as before:
\beqy
4\pi\log\int_\Si h_1e^{\chi_kw_1}\mrm dV_g&\le&4\pi\int_\Si\chi_kw_1\mrm dV_g+\int_\Si\fr{1}4|\n(\chi_kw_1)|^2\mrm dV_g+C\\
&\le&4\pi\|\chi_k\|_{L^\infty(\Si)}\int_\Si|w_1|\mrm dV_g+\int_\Si Q_{B_2}(\chi_kw)\mrm dV_g\\
&\le&(1+\e)\int_{\O'_k} Q_{B_2}(w)\mrm dV_g+\e\int_\Si Q_{B_2}(u)\mrm dV_g+C\\
&\le&\int_{\O'_k} Q_{B_2}(w)\mrm dV_g+\e\int_\Si Q_{B_2}(u)\mrm dV_g+C,
\eeqy
therefore
\beqa
\nonumber4\pi\l(\log\int_\Si e^{u_1}\mrm dV_g-\int_\Si u_1\mrm dV_g\r)&\le&4\pi\log\int_\Si e^{\chi_kw_1}\mrm dV_g+4\pi\|v_1\|_{L^\infty(\Si)}+4\pi\log\fr{1}{\d'}\\
\nonumber&\le&\int_{\O'_k} Q_{B_2}(w)\mrm dV_g+\e\int_\Si Q_{B_2}(u)\mrm dV_g+\e\int_\Si Q_{B_2}(u)\mrm dV_g+C\\
\label{u1}&\le&\int_{\O'_k} Q_{B_2}(w)\mrm dV_g+\e'\int_\Si Q_{B_2}(u)\mrm dV_g+C.
\eeqa
Putting together \eqref{u1u2} and \eqref{u1} and exploiting the fact that $\O'_k\cap\O'_{k'}=\es$ for any $k\ne k'$ we obtain
\beqy
&&4\pi K_1\l(\log\int_\Si e^{u_1}\mrm dV_g-\int_\Si u_1\mrm dV_g\r)+2\pi K_2\l(\log\int_\Si e^{u_2}\mrm dV_g-\int_\Si u_2\mrm dV_g\r)\\
&\le&\sum_{k=1}^{K_1}\l(\int_{\O'_k}Q_{B_2}(w)\mrm dV_g+\e\int_\Si Q_{B_2}(u)\mrm dV_g+C\r)\\
&\le&\int_\Si Q_{B_2}(w)\mrm dV_g+\e'\int_\Si Q_{B_2}(u)\mrm dV_g+C\\
&\le&\l(1+\e'\r)\int_\Si Q_{B_2}(w)\mrm dV_g+C,
\eeqy
which, up to re-naming $\e$, completes the proof.\\

The improved inequality concerning $\int_\Si Q_{G_2}(u)\mrm dV_g$ can be proved in the very same way.\\
This argument also works when $K_2>K_1$: to adapt it, just exchange the roles of $u_1$ and $u_2$ and write, in place of \eqref{xy}:
$$\fr{|x|^2}2+\fr{x\cd y}2+\fr{|y|^2}4=\fr{|y|^2}8+\fr{|2x+y|^2}8\q\q\q|x|^2+x\cd y+\fr{|y|^2}3=\fr{|y|^2}{12}+\fr{|2x+y|^2}4$$
\epf\

\section{Conclusion}\

The proof of Theorem \ref{exmolt} will follow easily by Lemma \ref{phipsi}, which enlightens the relation between low sub-levels of the energy functional and the join of barycenters, and by known results presented in Section $2$.\\

\blem\label{phipsi}${}$\\
Suppose $\rho\in(4K_1\pi,4(K_1+1)\pi)\x(4K_2\pi,4(K_2+1)\pi)$.\\
Then, for $L\gg0$ large enough there exist maps
$$\Phi_{B_2}:(\g_1)_{K_1}\star(\g_2)_{K_2}\to J_{B_2,\rho}^{-L}\q\q\q\Psi_{B_2}:J_{B_2,\rho}^{-L}\to(\g_1)_{K_1}\star(\g_2)_{K_2}$$
$$\Phi_{G_2}:(\g_1)_{K_1}\star(\g_2)_{K_2}\to J_{G_2,\rho}^{-L}\q\q\q\Psi_{G_2}:J_{G_2,\rho}^{-L}\to(\g_1)_{K_1}\star(\g_2)_{K_2}$$
such that $\Psi_{B_2}\c\Phi_{B_2}$ and $\Psi_{G_2}\c\Phi_{G_2}$ are homotopically equivalent to $\mrm{Id}_{(\g_1)_{K_1}\star(\g_2)_{K_2}}$.
\elem\

\bpf[Proof of Theorem \ref{exmolt}]${}$\\
For simplicity we will consider only the case of $B_2$.\\
From Lemma \ref{phipsi}, the map $\Phi_{B_2}$ induces an immersion of homology groups
$$H_q((\g_1)_{K_1}\star(\g_2)_{K_2})\sr{\l(\Phi_{B_2}\r)_{*,q}}\inc H_q\l(J_{B_2,\rho}^{-L}\r),$$ therefore, by Proposition \ref{hom},
$$b_{2K_1+2K_2-1}\l(J_{B_2,\rho}^{-L}\r)\ge b_{2K_1+2K_2-1}((\g_1)_{K_1}\star(\g_2)_{K_2})>0,$$
hence it is not contractible.\\
Suppose now that the system \eqref{b2toda} has no solutions. Then, by Lemma \ref{deform}, $J_{B_2,\rho}^{-L}$ would be a deformation retract of $J_{B_2,\rho}^L$ for any $L>0$. On the other hand, Corollary \ref{contr} implies that $J_{B_2,\rho}^L$ is contractible for $L$ large enough, so $J_{B_2,\rho}^{-L}$ would also be contractible. This is a contradiction.\\

To get a generic multiplicity result, we first use Lemma \ref{dense} to be able to use Morse theory for a generic choice of initial data; then, we just apply Lemma \ref{morsein} and Proposition \ref{hom}:
$$\#\tx{ solutions of }\eqref{b2toda}\ge\sum_{q=0}^{+\infty}\wt b_q\l(J_{B_2,\rho}^{-L}\r)\ge\sum_{q=0}^{+\infty}\wt b_q((\g_1)_{K_1}\star(\g_2)_{K_2})\ge\bin{K_1+\l\lfloor\fr{-\chi(\Si)}2\r\rfloor}{\l\lfloor\fr{-\chi(\Si)}2\r\rfloor}\bin{K_2+\l\lfloor\fr{-\chi(\Si)}2\r\rfloor}{\l\lfloor\fr{-\chi(\Si)}2\r\rfloor}$$
\epf\

To prove Lemma \ref{phipsi} we need a technical estimate concerning the test functions introduced in Theorem \ref{testfun}.

\blem\label{dist}${}$\\
Suppose $\rho\in(4K_1\pi,4(K_1+1)\pi)\x(4K_2\pi,4(K_2+1)\pi)$.\\
Then, for any $\z\in(\g_1)_{K_1}\star(\g_2)_{K_2}$ one has
$$\fr{1}{C\max\{1,\la(1-t)\}}\le d\l(f_{1,\Phi_{B_2}^\la(\z)},(\Si)_{K_1}\r)\le\fr{C}{\max\{1,\la(1-t)\}},$$
$$\fr{1}{C\max\{1,\la(1-t)\}}\le d\l(f_{1,\Phi_{G_2}^\la(\z)},(\Si)_{K_1}\r)\le\fr{C}{\max\{1,\la(1-t)\}},$$
$$\fr{1}{C\max\{1,\la t\}}\le d\l(f_{2,\Phi_{B_2}^\la(\z)},(\Si)_{K_2}\r)\le\fr{C}{\max\{1,\la t\}},$$
$$\fr{1}{C\max\{1,\la t\}}\le d\l(f_{2,\Phi_{G_2}^\la(\z)},(\Si)_{K_2}\r)\le\fr{C}{\max\{1,\la t\}}.$$
Moreover, if $t<1$ one has
$$f_{1,\Phi_{B_2}^\la(\z)}\us{\la\to+\infty}\wk\s_1,$$
whereas if $t>0$ one has
$$f_{2,\Phi_{B_2}^\la(\z)}\us{\la\to+\infty}\wk\s_2.$$
\elem\

\bpf${}$\\
We will prove only the first of the four former estimates, since the arguments needed are essentially the same. Clearly, such a proof is trivial when $\la(1-t)\le1$.\\
To get the upper estimate, we will show that $d\l(f_{1,\Phi_{B_2}^\la},\s_1^\la\r)\le\fr{C}{\la(1-t)}$, where
$$\s_1^\la:=\sum_{k=1}^{K_1}t_{1k}^\la\d_{x_{1k}}\q\q\q t_{1k}^\la=t_{1k}\fr{\int_\Si\fr{h_1}{\l(1+\l(\la(1-t)d\l(\cd,x_{1k}\r)\r)^2\r)^2}e^{-\fr{\ph_2}2}\mrm dV_g}{\int_\Si h_1e^{\ph_1-\fr{\ph_2}2}\mrm dV_g}.$$
From Lemma \ref{exp} we get, for any $\phi\in\mrm{Lip}(\Si)$ with $\|\phi\|_{\mrm{Lip}(\Si)}\le1$:
\beqy
\l|\int_\Si\l(f_{1,\Phi_{B_2}^\la}-\s_1^\la\r)\phi\mrm dV_g\r|&=&\fr{\int_\Si\l(h_1e^{\ph_1-\fr{\ph_2}2}-\int_\Si h_1e^{\ph_1-\fr{\ph_2}2}\mrm dV_g\s_1^\la\r)\phi\mrm dV_g}{\int_\Si h_1e^{\ph_1-\fr{\ph_2}2}\mrm dV_g}\\
&\le&C\fr{(\la(1-t))^2}{\max\{1,\la t\}^2}\int_\Si h_1e^{\ph_1-\fr{\ph_2}2}\l|\phi-\sum_k{t_{1k}^\la}\phi(x_{1k})\r|\mrm dV_g\\
&\le&C\fr{(\la(1-t))^2}{\max\{1,\la t\}^2}\sum_kt_{1k}\int_\Si\fr{\fr{h_1}{\l(1+\l(\la(1-t)d\l(\cd,x_{1k}\r)\r)^2\r)^2}}{\sum_{k'}\fr{t_{2k'}}{1+\l(\la td\l(\cd,x_{2k'}\r)\r)^2}}\mrm dV_g|\phi-\phi(x_{1k})|\mrm dV_g\\
&\le&C(\la(1-t))^2\sum_kt_{1k}\int_\Si\fr{d(\cd,x_{1k})}{\l(1+\l(\la(1-t)d\l(\cd,x_{1k}\r)\r)^2\r)^2}\mrm dV_g.
\eeqy
To conclude the estimate we need to show
$$\int_\Si\fr{d(\cd,x_{1k})}{\l(1+\l(\la(1-t)d\l(\cd,x_{1k}\r)\r)^2\r)^2}\mrm dV_g\le\fr{C}{(\la(1-t))^3}.$$
As in the proof of Lemma \ref{exp}, we easily find
$$\int_{\Si\sm B_\d(x_{1k})}\fr{d(\cd,x_{1k})}{\l(1+\l(\la(1-t)d\l(\cd,x_{1k}\r)\r)^2\r)^2}\mrm dV_g\le\fr{C}{(\la(1-t))^4}.$$
For the integral inside the disk, we use normal coordinates and a change of variables:
$$\int_\Si\fr{d(\cd,x_{1k})}{\l(1+\l(\la(1-t)d\l(\cd,x_{1k}\r)\r)^2\r)^2}\mrm dV_g\le\fr{C}{(\la(1-t))^3}\int_{B_{\la(1-t)\d}}\fr{|y|}{\l(1+|y|^2\r)^2}\mrm dy\le\fr{C}{(\la(1-t))^3}.$$\

To have a lower bound, we suffice to prove that, for any $\s=\s^\la\in(\Si)_{K_1}$ there exists $\phi_\s\in\mrm{Lip}(\Si)$ with $\|\phi_\s\|_{\mrm{Lip}(\Si)}\le1$ and
$$\l|\int_\Si\l(f_{1,\Phi_{B_2}^\la}-\s\r)\phi\mrm dV_g\r|\ge\fr{1}{C\max\{1,\la(1-t)\}}.$$
Precisely, by choosing
$$\phi_\s:=\min_{k'}d(\cd,x_{k'})\q\q\q\tx{if }\s=\sum_{k'=1}^{K_1}t_{k'}\d_{x_{k'}},$$
we obtain
\beqy
\l|\int_\Si\l(f_{1,\Phi_{B_2}^\la}-\s\r)\phi_\s\mrm dV_g\r|&=&\fr{\int_\Si\l(h_1e^{\ph_1-\fr{\ph_2}2}-\int_\Si h_1e^{\ph_1-\fr{\ph_2}2}\mrm dV_g\s\r)\phi_\s\mrm dV_g}{\int_\Si h_1e^{\ph_1-\fr{\ph_2}2}\mrm dV_g}\\
&=&\fr{\int_\Si h_1e^{\ph_1-\fr{\ph_2}2}\min_{k'}d(\cd,x_{k'})\mrm dV_g}{\int_\Si h_1e^{\ph_1-\fr{\ph_2}2}\mrm dV_g}\\
&\ge&\fr{(\la(1-t))^2}{C\max\{1,\la t\}^2}\int_\Si h_1e^{\ph_1-\fr{\ph_2}2}\min_{k'}d(\cd,x_{k'})\mrm dV_g\\
&\ge&\fr{(\la(1-t))^2}C\sum_{k=1}^{K_1}t_{1k}\int_{B_\d(x_{1k})}\fr{\min_{k'}d(\cd,x_{k'})}{\l(1+\l(\la(1-t)d\l(\cd,x_{1k}\r)\r)^2\r)^2}\mrm dV_g\\
&\ge&\fr{1}{C\la(1-t)}\int_{B_{\la(1-t)\d}(0)}\fr{\min_k\l|y-\la(1-t)x_{k'}\r|}{\l(1+|y|^2\r)^2}\mrm dy.
\eeqy
To conclude the proof, we just suffice to show that the last integral is uniformly bounded from below. The minimum will be attained by $x'=x_{k'}$ on a portion of the ball which measures at least $\fr{1}K$ of the whole ball.\\
By choosing $x'=x'^\la$ satisfying $\la x_1^\la\us{\la\to+\infty}\lto+\infty$, the integral goes to $+\infty$; otherwise, as shown in the proof of the first part of the Lemma, it is uniformly bounded.\\
To get the last claim, we need to show that $t_{ik}^\la\us{\la\to0}\to t_{ik}$, which in turn will follow from
$$\fr{\int_\Si\fr{h_1}{\l(1+\l(\la(1-t)d\l(\cd,x_{1k}\r)\r)^2\r)^2}e^{-\fr{\ph_2}2}\mrm dV_g}{\int_\Si\fr{h_1}{\l(1+\l(\la(1-t)d\l(\cd,x_{1k'}\r)\r)^2\r)^2}e^{-\fr{\ph_2}2}\mrm dV_g}\us{\la\to0}\to1\q\q\q\fa\,k,k'=1,\ds,K_1.$$
To get the last claim, we use the fact (proved in Lemma \ref{exp}) that the two integrals attains most of their mass around $x_{1k}$ and $x_{1k'}$, respectively: for any fixed $\d>0$ we have
\beqy\fr{\int_\Si\fr{h_1}{\l(1+\l(\la(1-t)d\l(\cd,x_{1k}\r)\r)^2\r)^2}e^{-\fr{\ph_2}2}\mrm dV_g}{\int_\Si\fr{h_1}{\l(1+\l(\la(1-t)d\l(\cd,x_{1k'}\r)\r)^2\r)^2}e^{-\fr{\ph_2}2}\mrm dV_g}&=&\fr{\int_{B_\d(x_{1k})}\fr{h_1}{\l(1+\l(\la(1-t)d\l(\cd,x_{1k}\r)\r)^2\r)^2}e^{-\fr{\ph_2}2}\mrm dV_g}{\int_{B_\d(x_{1k'})}\fr{h_1}{\l(1+\l(\la(1-t)d\l(\cd,x_{1k'}\r)\r)^2\r)^2}e^{-\fr{\ph_2}2}\mrm dV_g}+o_\la(1)\\
&=&\fr{\int_{B_\d(x_{1k})}\fr{\mrm dV_g}{\l(1+\l(\la(1-t)d\l(\cd,x_{1k}\r)\r)^2\r)^2}}{\int_{B_\d(x_{1k'})}\fr{\mrm dV_g}{\l(1+\l(\la(1-t)d\l(\cd,x_{1k'}\r)\r)^2\r)^2}}+o_\la(1).
\eeqy
By taking $\d$ close enough to $1$, each integral appearing in the last formula will be close to a Euclidean one, which will not depend on its center; therefore, by taking $\d_\la$ which goes to $0$ slow enough the ratio will tend to 1.
\epf\

\bpf[Proof of Lemma \ref{phipsi}]${}$\\
We will show the proof for $\Phi_{B_2}$ and $\Psi_{B_2}$.\\
Take $C$ as in Lemma \ref{dist}, $\e_0$ as in Lemma \ref{psik} and apply Theorem \ref{close} with $\e:=\fr{\e_0}{C^2}$. Then take $L:=L(\e)$ as in Theorem \ref{close} and define $\Phi_{B_2}:=\Phi_{B_2}^{\la_0}$ with $\Phi_{B_2}^{\la_0}((\g_1)_{K_1}\star(\g_2)_{K_2})\sub J_{B_2,\rho}^{-L}$.\\
To build the map $\Psi_{B_2}$ we start by the parameter $t'$:
$$t'(u)=t'(d_1(u),d_2(u)):=\l\{\bll0&\tx{if }d_2(u)\ge\e\\\fr{\e-d_2}{2\e-d_1-d_2}&\tx{if }d_1(u),d_2(u)<\e\\1&\tx{if }d_1(u)\ge\e\earr\r.\q\q\q\tx{with }d_i(u)=d_{\mrm{Lip}'}\l(f_{i,u},(\Si)_{K_i}\r).$$
$t'$ is well-defined because, by Theorem \ref{dist}, either $d_1(u)$ or $d_2(u)$ is less than $\e$.\\
Consider now $\psi_{K_1},\psi_{K_2}$ as in Lemma \ref{psik} and the push-forward $(\Pi_i)_*$ of the retractions $\Pi_i:\Si\to\g_i$, and define
$$\Psi_{B_2}(u):=((\Pi_1)_*(\psi_{K_1}(f_{1,u})),((\Pi_2)_*\psi_{K_2}(f_{2,u})),t'(u)).$$
This map is well-defined because if $\psi_{K_1}(f_{1,u})$ cannot be defined, then $d_1(u)\ge\e_0>\e$, hence $t'=1$. Similarly, we have $t'=0$ if $\psi_{K_2}(f_{2,u})$ is not defined.\\

Let us now construct the homotopical equivalence.\\
From the last assertion of Lemma \ref{phipsi}, $f_{i,\Phi_{G_2}^\la(\z)}$ converges to $\s_i$, and the convergence is preserved by the retractions $\psi_i$ and $(\Pi_i)_*$. However, the parameter $t'$ could be different from $t$.\\
Because of this, the homotopy map will consist in two steps: first we let $\la$ go to $+\infty$, then we pass from $t'$ to $t$.\\
Writing $\Psi_{B_2}\l(\Phi_{B_2}^\la(\z)\r)=\l((\Pi_1)_*\psi_{K_1}^\la(\z),(\Pi_2)_*\psi_{K_2}^\la(\z),t'^\la(\z)\r)$, the map $F:(\g_1)_{K_1}\star(\g_2)_{K_2}\x[0,1]\to(\g_1)_{K_1}\star(\g_2)_{K_2}$ will be given by $F:=F_2\ast F_1$, where
\beqy
F_1:(\z,s)=((\s_1,\s_2,t),s)&\mapsto&\l((\Pi_1)_*\l(\psi_{K_1}^\fr{\la_0}{1-s}(\z)\r),(\Pi_2)_*\l(\psi_{K_2}^\fr{\la_0}{1-s}(\z)\r),t'^{\la_0}(\z)\r)\\
F_2:\l(\l(\s_1,\s_2,t'^\la(\z)\r),s\r)&\mapsto&\l(\s_1,\s_2,(1-s)t'^{\la_0}(\z)+st\r).
\eeqy
Let us verify that $F_1$ is well defined: from Lemma \ref{dist}, $\psi_{K_1}$ is defined as long as $d_1\l(\Phi_{B_2}^\fr{\la_0}{1-s}(\z)\r)\le\e_0$, but if this does not occur, then
$$d_1\l(\Phi_{B_2}^{\la_0}(\z)\r)\ge\fr{1}{C\max\{1,\la_0(1-t)\}}\ge\fr{1}{C\max\l\{1,\fr{\la_0}{1-s}(1-t)\r\}}\ge\fr{d_1\l(\Phi_{B_2}^\fr{\la_0}{1-s}(\z)\r)}{C^2}\ge\e,$$
therefore $t'=1$; similarly, if $\psi_{K_2}$ is not defined then $t'=0$, hence $F_1$ makes sense.\\
Concerning $F_2$, if the first element in the join is not defined, then $t=1$ but one also gets $t'=1$, so there are no issues in their convex combination; similarly, if the second element is not defined, then $t=t'=0$, hence everything still works.
\epf\

\section{Extension to general systems}\

In Section $6$ of this paper we will give a proof of Theorem \ref{sublev} and of its Corollary \ref{exist} based on the same arguments of Sections $3,4,5$.\\
Since all the proofs are quite similar to the ones in the previous sections, we will be sketchy.\\

As a first thing, we consider Lemma \ref{retra}.\\
It is easy to see that, in the case of an $N$-dimensional system, we cannot take $N$ disjointed wedge sums of $g$ circles each: in fact, if we take $\g_1,\g_2$ in this way, then $\Si\sm(\g_1\cup\g_2)$ is a disjointed union of $g+1$ annuli, which cannot be retracted on a wedge sum of $g$ circles for $g>1$.\\
Anyway, to show the non-contractibility of low sub-levels, we suffice to take $\g_i$ as simple closed curves. In this way, we can easily build an arbitrary number of retractions on such disjointed curves; moreover, up to small perturbations, this can be done in such a way that all circles do not contain any of the points $p_m$.
The counterpart of Lemma \ref{retra} we need is the following:\\

\blem\label{retragen}${}$\\
Let $\Si$ be a surface with $\chi(\Si)\le0$.\\
Then, for any $N\in\N$ and for any given $p_1,\ds,p_M\in\Si$, there exist disjointed simple closed curves $\g_1,\ds,\g_N$ such that $p_m\nin\g_i$ for any $i,m$, and global retractions $\Pi_i:\Si\to\g_i$ for $i=1,\ds,N$.
\elem\

We will then need a sort of \qm{iterated join} of all the barycenters $(\g_i)_{K_i}$.\\
Actually, rather than simply repeating $N-1$ times the construction in Section $2$, we can equivalently consider the space $\bigstar_{i=1}^NX_i$ defined by
$$\bigstar_{i=1}^NX_i:=\fr{\Pi_{i=1}^NX_i\x\D^N}\sim,$$
where $\D_N$ is the unit simplex and $\sim$ is given by
$$(x_1,\ds,x_i,\ds,x_N,t_1,\ds,\ub{0}_{i},\ds,t_N)\sim(x_1,\ds,x'_i,\ds,x_N,t_1,\ds,\ub{0}_{i},\ds,t_N),$$
$$\fa\,x_1\in X_1,\,\ds,\,x_i,x'_i\in X_i,\,\ds,\,x_N\in X_N,\,(t_1,\ds,t_{i-1},t_{i+1},\ds,t_N)\in\D^N.$$
Such a space can be easily verified to be homeomorphic to $(\ds((X_1\star X_2)\star X_3)\star\ds\star X_N)$: writing $(x_1,\ds,x_N,t_1,\ds,t_N)\in\bigstar_{i=1}^NX_i$ and, for $i=2,\ds,N$, $\l(x'_i,x'_{i+1},t'_{i+1}\r)\in X_i\star X_{i+1}$, such an equivalence can be obtained by
$$\l\{\bll x'_i=x_i&i=1,\ds,N\\t'_i=\fr{t_i}{\sum_{j=1}^it_j}&i=2,\ds,N\earr\r.\q\q\q\iff\q\q\q\l\{\bll x_i=x'_i&i=1,\ds,N\\t_i=t'_i\Pi_{j=i+1}^N\l(1-t'_j\r)&i=2,\ds,N\earr\r..$$
Therefore, we will be considering the space
$$\mcal X:=\bigstar_{i=1}^N(\g_i)_{K_i}.$$
By the equivalence we just showed and Proposition \ref{hom}, we can summarize the properties of $\mcal X$ in the following lemma:\\

\blem\label{iterjoin}${}$\\
Let $\g_1,\ds,\g_N$ simple closed curves.\\
Then,
$$H_q(\mcal X)=\l\{\bll\Z&\tx{if }q=0,2\sum_{i=1}^NK_i-1\\0&\tx{if }q\ne0,2\sum_{i=1}^NK_i-1\earr\r.$$
In particular, it is not contractible.
\elem\

Concerning Moser-Trudinger inequality, we will use the general form of Theorem \ref{mt}, in the version originally proved in \cite{bat2}.\\

\bthm\label{mtgen}${}$\\
There exists $C>0$ such that for any $u\in H^1(\Si)^N$ one has
$$\sum_{i=1}^N\fr{8\pi}{a_{ii}}\l(\log\int_\Si e^{u_i}\mrm dV_g-\int_\Si u_i\mrm dV_g\r)\le\int_\Si Q_A(u)\mrm dV_g+C$$
In particular, \eqref{ja} is bounded from below if and only if $\rho_i\le\fr{8\pi}{a_{ii}}$ for all $i=1,\ds,N$.\\
Moreover, it coercive if and only if $\rho_i<\fr{8\pi}{a_{ii}}$ for all $i=1,\ds,N$. In this case, \eqref{general} has a minimizing solution.
\ethm\

The rest of this Section will be divided in three sub-sections, each of which is devoted to adapt the argument from Section $3$, $4$, $5$, respectively.\\

\subsection{Test functions}\

The test function we will consider are very similar to the ones in Theorem \ref{testfun}. We will still consider linear combinations of the \emph{bubbles} $\ph_i$, with coefficients depending on the entries of $A$.\\

\bthm\label{testgen}${}$\\
Define, for any $\la>0$ and $\z=(\s_1,\ds,\s_N,t_1,\ds,t_N)=\l(\sum_{k=1}^{K_1}t_{1k}\d_{x_{1k}},\ds,\sum_{k=1}^{K_N}t_{Nk}\d_{x_{Nk}},t_1,\ds,t_N\r)\in\mcal X,$
$$\ph_i=\ph_i^\la(\z)=\log\sum_{k=1}^{K_i}\fr{t_{ik}}{\l(1+(\la t_id(\cd,x_{ik}))^2\r)^2}\q i=1,\ds,N;\q\q\q\Phi_A^\la(\z)=\l(\sum_{i=1}^N\fr{a_{1i}}{a_{ii}}\ph_i,\ds,\sum_{i=1}^N\fr{a_{Ni}}{a_{ii}}\ph_i\r).$$
If $\rho_i\in\l(\fr{8\pi}{a_{ii}}K_i\pi,\fr{8\pi}{a_{ii}}(K_i+1)\pi\r)$ for $i=1,\ds,N$, then
$$J_{A,\rho}\l(\Phi_A^\la(\z)\r)\us{\la\to+\infty}\to-\infty\q\q\q\tx{uniformly in }\z\in\mcal X.$$
\ethm\

\bpf${}$\\
We start by estimate the term involving $Q_A$: we write
$$Q_A\l(\Phi_A^\la(\z)\r)=\fr{1}2\sum_{i=1}^N\fr{1}{a_{ii}}|\n\ph_i|^2+\sum_{i<j}\fr{a_{ij}}{a_{ii}a_{jj}}\n\ph_i\cd\n\ph_j$$
and exploit that, as in \eqref{grad1}, \eqref{grad12}, \eqref{grad2},
$$\int_\Si|\n\ph_i|^2\mrm dV_g\le32\pi\fr{K_i}{a_{ii}}\log\max\{1,\la t_i\}+C\q\q\q\int_\Si\n\ph_i\cd\n\ph_j\mrm dV_g=O(1);$$
therefore
\bequ\label{q}
\int_\Si Q_A\l(\Phi_A^\la(\z)\r)\mrm dV_g\le16\pi\sum_{i=1}^N\fr{K_i}{a_{ii}}\log\max\{1,\la t_i\}+C.
\eequ
Concerning averages, the very same argument of Lemma \ref{average} yields
\bequ\label{av}
\int_\Si\ph_i\mrm dV_g=-4\log\max\{1,\la t_i\}+O(1).
\eequ
Finally, since $d(x_{ik},p_m)\ge\d>0$, similar computations as in Lemma \ref{exp} show that
\beqy
\int_\Si\wt h_ie^{\sum_{j=1}^N\fr{a_{ij}}{a_{jj}}\ph_j}\mrm dV_g&\sim&\prod_{j\ne i}\max\{1,\la t_i\}^{-4\fr{a_{ij}}{a_{jj}}}\sum_{k=1}^{K_i}t_{ik}\int_{B_\d(x_{ik})}\fr{\mrm dV_g}{\l(1+(\la t_id(\cd,x_{ik}))^2\r)^2}\\
&\sim&\max\{1,\la t_j\}^2\prod_{j\ne i}\max\{1,\la t_j\}^{-4\fr{a_{ij}}{a_{jj}}},
\eeqy
hence
\bequ\label{log}
\log\int_\Si\wt h_ie^{\sum_{j=1}^N\fr{a_{ij}}{a_{jj}}\ph_j}\mrm dV_g=-2\log\max\{1,\la t_i\}-4\sum_{j\ne i}\fr{a_{ij}}{a_{jj}}\log\max\{1,\la t_j\}+O(1).
\eequ
The Lemma now follows from \eqref{q}, \eqref{av}, \eqref{log}.
\epf\

\subsection{Improved M-T inequalities}\

This sub-section is devoted to characterize low sub-levels of $J_{A,\rho}$ in the same way as Theorem \ref{close}.\\
Namely, we will show that, as $J_{A,\rho}(u)\ll0$, at least one of its components is close to the corresponding barycenters space.\\

\bthm\label{closegen}${}$\\
Suppose $\rho_i\in\l(\fr{8\pi}{a_{ii}}K_i\pi,\fr{8\pi}{a_{ii}}(K_i+1)\pi\r)$ for $i=1,\ds,N$.\\
Then, for any $\e>0$ there exists $L=L_\e>0$ such that any $u\in J_{A,\rho}^{-L}$ verifies, for some $i=1,\ds,N$,
$$d_{\mrm{Lip}'}(f_{i,u},(\Si)_{K_i})\le\e.$$
\ethm\

As usual, such a result will follow, via standard arguments, from an improved Moser-Trudinger inequality, which has the following form:\\

\blem\label{imprgen}${}$\\
Let $\d>0$, $K_1,\ds,K_N\in\N$, $\{\O_{ik}\}_{i=1,\ds,N,\,k=1,\ds,K_i}$ satisfy
\beqy
d(\O_{ik},\O_{ik'})\ge\d&\q\q\q&\fa\,i=1,\ds,N,\,k,k'=1,\ds,K_i,\,k\ne k'\\
\int_{\O_{ik}}f_{i,u}\mrm dV_g\ge\d&\q\q\q&\fa\,i=1,\ds,N,\,k=1,\ds,K_i.\\
\eeqy
Then, for any $\e>0$ there exists $C>0$, not depending on $u$, such that
\bequ\label{impr}
\sum_{i=1}^N\fr{8\pi}{a_{ii}}K_i\l(\log\int_\Si e^{u_i}\mrm dV_g-\int_\Si u_i\mrm dV_g\r)\le(1+\e)\int_\Si Q_A(u)\mrm dV_g+C.
\eequ
\elem\

To prove Lemma \ref{imprgen}, we will need a couple of ingredients.\\
First of all, a covering Lemma slightly different from Lemma \ref{cover}.\\

\blem\label{covergen}${}$\\
Let $\d>0$, $K_1,\ds,K_N\in\N$ be given, $f_1,\ds,f_N\in L^1(\Si)$ be positive a.e. and such that $\|f_1\|_{L^1(\Si)}=\ds=\|f_N\|_{L^1(\Si)}=1$ and $\{\O_{ik}\}_{i=1,\ds,N,k=1,\ds,K_i}$ satisfy
\beqy
d(\O_{ik},\O_{ik'})\ge\d&\q\q\q&\fa\,i=1,\ds,N,\,k,k'=1,\ds,K_i,\,k\ne k'\\
\int_{\O_{ik}}f_{i,u}\mrm dV_g\ge\d&\q\q\q&\fa\,i=1,\ds,N,\,k=1,\ds,K_i.\\
\eeqy
Then, there exists $\d'=\d'(\d,N,\Si)$, $K\in\l\{1,\ds,\sum_{i=1}^NK_i\r\}$, $\mcal K_i\sub\{1,\ds,K\}$ with $|\mcal K_i|=K_i$ and $\{\O_k\}_{k=1}^K$ such that
\beqy
d(\O_k,\O_{k'})\ge\d'&\q\q\q&\fa\,k,k'=1,\ds,K,\,k\ne k',\\
\int_{\O_k}f_i\mrm dV_g\ge\d'&\q\q\q&\fa\,i=1,\ds,N,\,k\in\mcal K_i.
\eeqy
\elem\

\bpf${}$\\
Fix $\d_0:=\fr{\d}{3N+2}$ and write, by compactness, $\Si=\Cup_{l=1}^LB_{\d_0}(x_l)$ forsome $x_1,\ds,x_L\in\Si$ and $L=L(\d_0,\Si)$. Then take $\{x_{ik}\}_{i=1,\ds,N,k=1,\ds,K_i}\sub\{x_l\}_{l=1}^L$ such that
$$\int_{B_{\d_0}(x_{ik})}f_i\mrm dV_g=\max\l\{\int_{B_{\d_0}(x_l)}f_i\mrm dV_g:\,B_{\d_0}(x_l)\cap\O_{ik}\ne\es\r\}.$$
We define inductively $\O'_{ik}$ as suitable unions of such balls in the following way.\\
Start by setting $\O'_{1k}:=B_{\d_0}(x_{1k})$ for $k=1,\ds,K_1$. Then, since by construction $d(x_{ik},x_{ik'})\ge3N\d_0$ for any $k\ne k'$, each $x_{2k'}$ could be $3\d_0$-close to at most one $x_{1k}$; up to re-labeling indices, we can assume that $x_{2k'}$ does not have a $3\d_0$-close $x_{1k}$ point if and only if $k=K_1+1,\ds,K'_2$; we then define
$$\O'_{2k}=\l\{\bll B_{\d_0}(x_{1k})\cup B_{\d_0}(x_{2k'})&\tx{if }k\le K_1,\,d(x_{1k},x_{2k'})<3\d'\tx{ for some }k'\\B_{\d_0}(x_{1k})&\tx{if }k\le K_1,\,d(x_{1k},x_{2k'})\ge3\d'\tx{ for any }k'\\B_{\d_0}(x_{2k})&\tx{if }k=K_1+1,\ds,K'_2\earr\r..$$
We then iterate the construction by setting
$$\O'_{i+1,k}:=\l\{\bll\O_{ik}\cup B_{\d_0}(x_{i+1,k'})&\tx{if }k\le K'_i,\,d(\O_{ik},x_{i+1,k'})<2\d'\tx{ for some }k'\\\O_{ik}&\tx{if }k\le K'_i,\,d(\O_{ik},x_{i+1,k'})\ge2\d'\tx{ for any }k'\\B_{\d_0}(x_{i+1,k})&\tx{if }k=K'_i+1,\ds,K'_{i+1}\earr\r.,$$
which is allowed since $d(\O_{i-1,k},\O_{i-1,k'})\ge(3(N-i+1)+1)\d'$. Finally, we set $\O_k:=\O_{N,k}$.\\
By construction one gets $d(\O_k,\O_{k'})\ge\d_0$ and, setting $\mcal K_i:=\{k:\,B_{\d_0}(x_{ik'})\sub\O_k\tx{ for some }k'\}$,
we get
$$\int_{\O_k}f_i\mrm dV_g\ge\int_{B_{\d_0}(x_{ik'})}f_i\mrm dV_g\ge\fr{\d}L,$$
hence the conclusion follows by setting $\d':=\min\l\{\d_0,\fr{\d}L\r\}$.
\epf\

We will also use the following property of positive definite matrices, which will help in dealing with Moser-Trudinger inequalities involving different matrices:\\

\blem\label{matrix}${}$\\
Let $A=(a_{ij})_{i,j=1,\ds,N}$ be a symmetric, positive definite matrix, $\mcal I\sub\{1,\ds,N\}$ a given set of indices, $B:=A|_{\mcal I\x\mcal I}$ a submatrix of $A$ and $a^{ij},b^{ij}$ the entries of the inverse $A^{-1},B^{-1}$ of $A,B$, respectively.\\
Then, for any $x=(x_1,\ds,x_N)\in\R^N$,
$$\sum_{i,j\in\mcal I}b^{ij}x_ix_j\le\sum_{i,j=1}^Na^{ij}x_ix_j.$$
\elem\

Such a result is equivalent to showing that the matrix $C$ defined by $c_{ij}:=\l\{\bll a^{ij}-b^{ij}&\tx{if }i,j\in\mcal I\\a^{ij}&\tx{otherwise}\earr\r.$ is positive semi-definite. This follows from $A$ being positive definite and $AC$ being positive semi-definite, since $(AC)_{ij}=\l\{\bll1&\tx{if }i=j\nin\mcal I\\0&\tx{otherwise}\earr\r.$.\\

\bpf[Proof of Lemma \ref{imprgen}]${}$\\
We start by applying Lemma \ref{covergen} to get the sets $\{\O_k\}_{k=1}^K$ and we take, for each of them, a corresponding cutoff $\chi_k$, as in \eqref{chi}.\\
We write $u_i-\int_\Si u_i\mrm dV_g=v_i+w_i$ with $v_i\in L^\infty(\Si)$ chosen by truncation in Fourier decomposition. Then, we fix $k$ and apply a Moser-Trudinger inequality to $\chi_kw$, taking account only of the components $w_i$ for which $k\in\mcal K_i$. We get:
\beqy
&&\sum_{\{i:\,k\in\mcal K_i\}:=\mcal I}\fr{8\pi}{a_{ii}}\l(\log\int_\Si\wt h_ie^{u_i}\mrm dV_g-\int_\Si u_i\mrm dV_g\r)\\
&\le&\sum_{i\in\mcal I}\fr{8\pi}{a_{ii}}\log\int_\Si\wt h_ie^{\chi_kw_i}\mrm dV_g+\e\int_\Si Q_A(u)\mrm dV_g+C\\
&\le&(1+\e)\int_{B_\fr{\d'}2(\O_k)}Q_{A|_{\mcal I\x\mcal I}}(w)\mrm dV_g+\e\int_\Si Q_{A|_{\mcal I\x\mcal I}}(u)\mrm dV_g+\e\int_\Si Q_A(u)\mrm dV_g+C\\
&\le&\int_{B_\fr{\d'}2(\O_k)}Q_A(w)\mrm dV_g+2\e\int_\Si Q_A(u)\mrm dV_g+C,
\eeqy
where the last passage uses Lemma \ref{matrix}. Taking a sum over all $k$'s proves the Lemma, being $|\mcal K_i|=K_i$.
\epf\

\subsection{Conclusion}\

The estimates on test functions from Theorem \ref{testgen} and the result from Theorem \ref{closegen} allow to get our claim via standard argument.

\bpf[Proof of Theorem \ref{sublev} and Corollary \ref{exist}]${}$\\
To prove Theorem \ref{sublev} we just need, as in Lemma \ref{phipsi}, two maps
$$\Phi_A:\mcal X\to J_{A,\rho}^{-L}\q\q\q\Psi_A:J_{A,\rho}^{-L}\to\mcal X$$
such that $\Psi_A\circ\Phi_A$ is homotopically equivalent to $\mrm{Id}_{\mcal X}$.\\
As a first thing, arguing as in Lemma \ref{dist} shows that
$$\fr{1}{C\max\{1,\la t_i\}}\le d\l(f_{i,\Phi_A^\la(\z)},(\Si)_{K_i}\r)\le\fr{C}{\max\{1,\la t_i\}}$$
and, if $t_i>0$,
$$f_{i,\Phi_A^\la(\z)}\us{\la\to+\infty}\wk\s_i.$$
We then take $\e_0$ as in Lemma \ref{psik}, fix $\e:=\fr{\e_0}{C^2}$, $L:=L(\e)$ as in Theorem \ref{closegen} and $\la_0$ such that $\Phi_A^{\la_0}(\mcal X)\sub J_{A,\rho}^{-L}$. Then, we define
$$\Psi_A(u):=\l((\Pi_1)_*(\psi_{K_1}(f_{1,u}),\ds,\psi_{K_N}(f_{N,u}),t'_1(u),\ds,t'_N(u)\r),\q\q\q t'_i(u):=\fr{(\e-d_i)^+}{\sum_{j=1}^N(\e-d_j)^+},$$
where we recall $d_{\mrm{Lip}'}\l(f_{i,u},(\Si)_{K_i}\r)$.\\
Writing $\Psi_A(\Phi_A^\la(\z))=\l((\Pi_1)_*\psi_{K_1}^\la(\z),\ds,(\Pi_N)_*\psi_{K_N}^\la(\z),{t'_1}^\la(\z),\ds,{t'_N}^\la(\z)\r)$, a suitable homotopical equivalence $F:\mcal X\x[0,1]\to\mcal X$ will be given by $F_2*F_1$, where
\beqy
F_1:(\z,s)=((\s_1,\ds,\s_N,t_1,\ds,t_N),s)&\mapsto&\l((\Pi_1)_*\psi_{K_1}^{\fr{\la_0}{1-s}}(\z),\ds,(\Pi_N)_*\psi_{K_N}^{\fr{\la_0}{1-s}}(\z),{t'_1}^{\la_0}(\z),\ds,{t'_N}^{\la_0}(\z)\r)\\
F_2:\l(\l(\s_1,\ds,\s_N,{t'_1}^{\la_0}(\z),\ds,{t'_N}^{\la_0}(\z)\r),s\r)&\mapsto&\l(\s_1,\ds,\s_N,(1-s){t'_1}^{\la_0}(\z)+st_1,\ds,(1-s){t'_N}^{\la_0}(\z)+st_N\r)
\eeqy
Arguing as in the proof of Lemma \ref{phipsi} shows that $F$ is well-posed and has all the desired properties.\\
Finally, if the compactness condition stated in Theorem \ref{sublev} holds, then Corollary \ref{exist} follows because of Lemma \ref{deform} and Corollary \ref{contr}, since both hold true for general systems.
\epf\

\begin{remark}${}$\\
After the submission of this manuscript, a new result \cite{lz} concerning compactness of solutions of \eqref{b2toda} and \eqref{g2toda} was published. It basically improves Theorem \ref{quant} and Corollary \ref{comp} and extends them to the case of singular $B_2,G_2$ systems, namely
$$\l\{\bl-\D u_1=2\rho_1\l(\fr{h_1e^{u_1}}{\int_\Si h_1e^{u_1}\mrm dV_g}-1\r)-\rho_2\l(\fr{h_2e^{u_2}}{\int_\Si h_2e^{u_2}\mrm dV_g}-1\r)-\sum_{m=1}^M\a_{1m}(\d_{p_m}-1)\\-\D u_2=2\rho_2\l(\fr{h_2e^{u_2}}{\int_\Si h_2e^{u_2}\mrm dV_g}-1\r)-2\rho_1\l(\fr{h_1e^{u_1}}{\int_\Si h_1e^{u_1}\mrm dV_g}-1\r)-\sum_{m=1}^M\a_{2m}(\d_{p_m}-1)\earr\r.,$$
$$\l\{\bl-\D u_1=2\rho_1\l(\fr{h_1e^{u_1}}{\int_\Si h_1e^{u_1}\mrm dV_g}-1\r)-\rho_2\l(\fr{h_2e^{u_2}}{\int_\Si h_2e^{u_2}\mrm dV_g}-1\r)-\sum_{m=1}^M\a_{1m}(\d_{p_m}-1)\\-\D u_2=2\rho_2\l(\fr{h_2e^{u_2}}{\int_\Si h_2e^{u_2}\mrm dV_g}-1\r)-3\rho_1\l(\fr{h_1e^{u_1}}{\int_\Si h_1e^{u_1}\mrm dV_g}-1\r)-\sum_{m=1}^M\a_{2m}(\d_{p_m}-1)\earr\r..$$
with $\a_{im}>-1$ and $p_m\in\Si$ given.\\
In \cite{lwyz}, Theorems $1.3,1.7$, the authors show that the local blow-up masses $\s_i(x)$ must be integer multiples of $4\pi$ if $x\notin\{p_1,\ds,p_m\}$ and $\s_i(p_m)=4\pi(N_{i,0}+N_{i,1}(1+\a_{1m})+N_{i,2}(1+\a_{2m}))$ for some integers $N_{ij}$. This implies that in Theorem \ref{conccomp} \emph{Compactness} must occur if neither $\rho_1$ nor $\rho_2$ belong to the the following discrete set $\La$:
$$\La:=4\pi\l\{N_0+\sum_{m=1}^MN_{1m}(1+\a_{1,m})+\sum_{m=1}^MN_{2m}(1+\a_{2m}),\,N_0,N_{1m},N_{2m}\in\N\r\}.$$
This fact, in view of Theorem \ref{sublev} and Corollary \ref{exist}, permits to extend Theorem \ref{exmolt} to the singular case (with coefficients $\a_{im}\ge0$) under the assumption $\rho_1,\rho_2\nin\La$. In case of no singularities, we get $\La=4\pi\N$, which means that Theorem \ref{exmolt} can be extended, for the regular $G_2$ Toda system, to any $\rho\nin4\pi\N\x\R\cup\R\x4\pi\N$ under no upper bound on $\rho_1,\rho_2$.\\
The results in \cite{lwyz} also hold true for the singular $A_2$ Toda system, thus giving an extension of some of the variational existence results from \cite{bjmr,bat1,batmal}
\end{remark}\

\section*{Acknowledgments}\

The authors wishes to thank to Department of Mathematics of University of British Columbia for the hospitality provided during his visit.\\
Particular gratitude is expressed to Prof. Juncheng Wei and Wen Yang for the discussions concerning the topics of the paper.\\

\bibliography{b2g2}
\bibliographystyle{abbrv}

\end{document}